\documentclass{amsart}

\usepackage{amsthm, amsmath, amssymb, url}

\newtheorem{thm}{Theorem}[section]

\newtheorem{prop}[thm]{Proposition}
\newtheorem{lemma}[thm]{Lemma}

\newtheorem{cor}[thm]{Corollary}

\theoremstyle{definition}
\newtheorem{alg}[thm]{Algorithm}
\newtheorem{rem}[thm]{Remark}
\newtheorem{ex}[thm]{Example}

\newcommand{\F}{\mathbb{F}}
\newcommand{\Q}{\mathbb{Q}}

\newcommand{\Z}{\mathbb{Z}}
\newcommand{\C}{\mathbb{C}}
\renewcommand{\O}{\mathcal{O}}
\newcommand{\pibar}{\overline{\pi}}
\newcommand{\p}{\mathfrak{p}}
\newcommand{\q}{\mathfrak{q}}
\newcommand{\rr}{\mathfrak{r}}
\newcommand{\B}{\mathcal{B}}

\newcommand{\ceil}[1]{\lceil #1 \rceil}
\newcommand{\abs}[1]{\left| #1 \right|}

\numberwithin{equation}{section}

\DeclareMathOperator{\Jac}{Jac}

\DeclareMathOperator{\Frob}{Frob}

\DeclareMathOperator{\End}{End}
\DeclareMathOperator{\Aut}{Aut}
\DeclareMathOperator{\Gal}{Gal}
\DeclareMathOperator{\lcm}{lcm}
\DeclareMathOperator{\Disc}{\Delta}

\DeclareMathOperator{\Norm}{Norm}

\title[Computing endomorphism rings of genus 2 Jacobians]{Computing endomorphism rings of Jacobians of genus 2 curves over finite fields}
\author{David Freeman}
\address{University of California, Berkeley \\
	{\tt dfreeman@math.berkeley.edu}}
\author{Kristin Lauter}
\address{Microsoft Research \\
	{\tt klauter@microsoft.com}}

\date{\today}

\begin{document}

\begin{abstract}
We present probabilistic algorithms which, given a genus $2$ curve $C$ defined over a finite field and a quartic CM field $K$, determine whether the endomorphism ring of the Jacobian $J$ of $C$ is the full ring of integers in $K$.  In particular, we present algorithms for computing the field of definition of, and the action of Frobenius on, the subgroups $J[\ell^d]$ for prime powers $\ell^d$.  We use these algorithms to create the first implementation of Eisentr\"ager and Lauter's algorithm for computing Igusa class polynomials via the Chinese Remainder Theorem \cite{el}, and we demonstrate the algorithm for a few small examples.  We observe that in practice the running time of the CRT algorithm is dominated not by the endomorphism ring computation but rather by the need to compute $p^3$ curves for many small primes $p$.
\end{abstract}

\maketitle

\section{Introduction}
Many public-key cryptographic protocols are based on the difficulty of the discrete logarithm problem in groups of points on elliptic curves and Jacobians of hyperelliptic curves.  For such protocols one needs to work in a subgroup of large prime order of the Jacobian of the curve, so it is useful to be able to construct curves over finite fields whose Jacobians have a specified number of points.  

The problem of constructing elliptic curves over finite fields with a given number of points has been studied extensively.  Current solutions rely on computing the $j$-invariant via the construction of the Hilbert class polynomial for a quadratic imaginary field.  There are three different approaches to computing the Hilbert class polynomial: a complex-analytic algorithm \cite{am}, \cite{enge}; a Chinese Remainder Theorem algorithm \cite{cnst}, \cite{alv}; and a $p$-adic algorithm \cite{couveignes-henocq}, \cite{broker}.  The best running time for these algorithms is $\tilde O(\abs{d})$, where $d$ is the discriminant of the quadratic imaginary field \cite{enge}, \cite{broker}.

Analogous methods exist for constructing genus 2 curves with a given number of points on their Jacobians.  In this case, the solutions rely on computing the curves' Igusa invariants via the computation of Igusa class polynomials for quartic CM fields.  Again there are three different approaches: a complex-analytic algorithm \cite{spallek}, \cite{vanwamelen}, \cite{weng}, \cite{cohn-lauter}; a Chinese Remainder Theorem algorithm \cite{el}; and a $p$-adic algorithm \cite{g2}.  These algorithms are less extensively developed than their elliptic curve analogues, and to date there is no running time analysis for any of them.

In this paper we study the implementation of Eisentr\"ager and Lauter's Chinese Remainder Theorem algorithm \cite{el}.  The algorithm takes as input a primitive quartic CM field $K$, i.e.~a purely imaginary quadratic extension of a real quadratic field with no proper imaginary quadratic subfields, and produces the Igusa class polynomials of $K$.  
The basic outline of the algorithm is as follows:
\begin{enumerate}
	\item \label{st:bound} Define $S$ to be a set of primes with certain splitting behavior in the field $K$ and its reflex field $K^*$. 
	\item For each prime $p$ in $S$:
	\begin{enumerate}
		\item For each triple $(i_1,i_2,i_3) \in \F_p^3$ of Igusa invariants, construct a genus $2$ curve $C$ over $\F_p$ corresponding to that triple.
		\item \label{st:end} Check the isogeny class of each curve. For each curve in the desired isogeny class, compute the endomorphism ring of the Jacobian of the curve and keep only those curves for which the endomorphism ring is the full ring of integers $\O_K$.
		\item Construct the Igusa class polynomials mod $p$ from the triples collected in Step \ref{st:end}.
	\end{enumerate}
	\item Use the Chinese Remainder Theorem or the Explicit CRT \cite{bernstein} to construct the Igusa polynomials either with rational coefficients or modulo a prime of cryptographic size. 
\end{enumerate}
One advantage of the CRT algorithm over other algorithms for computing Igusa class polynomials is that the CRT algorithm does not require that the real quadratic subfield have class number one.

Our contribution is to provide an efficient probabilistic algorithm for computing endomorphism rings of Jacobians of genus 2 curves over small prime fields.  Using this algorithm to compute endomorphism rings, we have implemented a probabilistic version of the full Eisentr\"ager-Lauter CRT algorithm (Algorithm \ref{a:crt}) in MAGMA and used it to compute Igusa class polynomials for several fields $K$ with small discriminant.

It was previously believed that computing endomorphism rings would be the bottleneck in the genus 2 CRT algorithm.  Our results are surprising in the sense that we find that the time taken to compute the endomorphism rings with our probabilistic algorithms is negligible compared with the time needed to compute $p^3$ genus 2 curves via Mestre's algorithm for each small prime $p$.  For example, for $K = \Q(i\sqrt{13 + 2\sqrt{13}})$ and $p=157$, the largest prime for which endomorphism rings are computed for this $K$, our (unoptimized) MAGMA program takes about 52 minutes to loop through $157^3$ curves and find $243$ curves in the specified isogeny class.  Our probabilistic algorithm (also implemented in MAGMA) applied to these $243$ curves then takes $16.5$ {\it seconds} to find the single curve whose Jacobian has endomorphism ring equal to $\O_K$.

The algorithm works as follows.  Let $C$ be a genus 2 curve over a finite field $\F_p$, and let $J$ be its Jacobian; we assume $J$ is ordinary.  
Let $K$ be a primitive quartic CM field, which we assume is given via an embedding in $\C$.
The first test is whether $\End(J)$, the endomorphism ring of $J$, is an order in $\O_K$.  This computation is outlined in \cite[Section 5]{el} and described in more detail in Section \ref{s:zeta} below.  If $\End(J)$ is an order in $\O_K$, we compute a set of possible elements $\pi \in \O_K$ that could represent the Frobenius endomorphism of $J$.  If $\pi$ represents the Frobenius endomorphism, then its complex conjugate $\pibar$ represents the Verschiebung endomorphism.  

We next determine a set $\{\alpha_i\}$ of elements of $\O_K$ such that $\Z[\pi,\pibar,\{\alpha_i\}] = \O_K$.  It follows that $\End(J) = \O_K$ if and only if each $\alpha_i$ is an endomorphism of $J$.  We show in Section \ref{s:generators} that we can take each $\alpha_i$ to have one of two forms: either $\alpha_i = \frac{\pi^k-1}{\ell}$ for some positive integer $k$ and prime $\ell$, or $\alpha_i = \frac{h_i(\pi)}{\ell^d}$ for some cubic polynomial $h_i$ with integer coefficients and some prime power $\ell^d$.  In Section \ref{s:fields} we show how to determine whether an element of the first form is an endomorphism; this is equivalent to determining the field of definition of the $\ell$-torsion points of $J$.  In Section \ref{s:frob} we show how to determine whether an element of the second form is an endomorphism; this is equivalent to computing the action of Frobenius on a basis of $J[\ell^d]$.  The main results are Algorithms \ref{a:fields} and \ref{a:frob}, two very efficient probabilistic algorithms which check fields of definition and compute the action of Frobenius, respectively.
The running times of these algorithms depend primarily on the sizes of the fields over which the points of $J[\ell^d]$ are defined.  Section \ref{s:time} provides upper bounds for these sizes in terms of the prime $\ell$ and the size of the base field $p$.

A detailed statement of the Eisentr\"ager-Lauter CRT algorithm, incorporating the algorithms of Sections \ref{s:zeta}, \ref{s:fields}, and \ref{s:frob}, appears in Section \ref{s:crt}.  Section \ref{s:notes} describes various ways in which we have modified our MAGMA implementation to improve the algorithm's performance.  Finally, in Section \ref{s:ex} we give examples of our algorithm run on several small quartic CM fields. 

\subsection*{Notation and assumptions} 

Throughout this paper, a {\it curve} will refer to a smooth, projective, absolutely irreducible algebraic curve $C$.  The Jacobian of $C$, denoted $\Jac(C)$, is an abelian variety of dimension $g$, where $g$ is the genus of $C$.  We assume throughout that $p$ is a prime, and that $\Jac(C)$ is an ordinary abelian variety modulo $p$.

A number field $K$ is a {\it CM field} if it is a totally imaginary quadratic extension of a totally real field.  We denote by $K^*$ the reflex field of $K$, and by $K_0$ the real quadratic subfield of $K$.
A CM field is {\it primitive} if it has no proper CM subfields. We will assume unless otherwise noted that $K$ is a primitive quartic CM field not isomorphic to $\Q(\zeta_5)$.  This implies that $K$ is either Galois cyclic or non-Galois.   If $K$ is Galois cyclic, then $K^*=K$; if $K$ is non-Galois, then $K^*$ is another primitive quartic CM field \cite[p. 64]{Shimura}.  A curve $C$ has {\it CM by $K$} if the endomorphism ring of $\Jac(C)$ is isomorphic to an order in $\O_K$, the ring of integers of the CM field $K$. 

\subsection*{Acknowledgments}

This research was conducted during the first author's internship at Microsoft Research, Redmond, during the summer of 2006.  The first author thanks Microsoft for its hospitality and Denis Charles, Jean-Marc Couveignes, and Edward Schaefer for many helpful discussions.  The second author thanks Pierrick Gaudry for helpful correspondence and pointers to his code.  Both authors thank Reinier Br\"oker, David Kohel, and Christophe Ritzenthaler for their feedback on previous versions of this paper.


\section{Computing zeta functions and the Frobenius element}
\label{s:zeta}

To determine whether the Jacobian $J$ of a given genus 2 curve $C$ has endomorphism ring equal to $\O_K$, the first step is to determine whether the endomorphism ring is even an order in $\O_K$.  This is accomplished by computing the characteristic polynomial of Frobenius, to see if the Frobenius element corresponds to an algebraic integer of $K$.  This in turn is equivalent to determining the zeta function of $C$, which can be computed by finding the number of points on the curve and its Jacobian, $n = \#C(\F_p)$ and $m = \#J(\F_p)$.  For a given field $K$ there are several possibilities for the pairs $(n, m)$, as described in \cite[Prop. 4]{el}.  

In this section we give an explicit algorithm that determines whether $\End(J)$ is an order in $\O_K$ and if so, gives a set $S \subset \O_K$ of possibilities for the Frobenius endomorphism of $J$.  The main point is to find the possible Frobenius elements by finding generators of certain principal ideals (Step \ref{st:gen}) with absolute value equal to $\sqrt{p}$ (Step \ref{st:unit}).

\begin{alg}
\label{a:zeta}
Let $K$ be a primitive quartic CM field and $K^*$ the reflex of $K$.
The following algorithm takes as input the field $K$, a prime $p$ that splits completely in $K$ and splits completely into principal ideals in $K^*$, and a curve $C$ defined over the finite field $\F_p$.  The algorithm returns {\tt true} or {\tt false} according to whether $\End(J)$ is an order in $\O_K$, where $J = \Jac(C)$.  If the answer is {\tt true}, the algorithm also outputs a set $S \subset \O_K$ that consists of the $\Aut(K/\Q)$-orbit of the Frobenius endomorphism of $J$.
\begin{enumerate}
\item Compute the decomposition $p = \p_1 \p_2 \p_3 \p_4$ in $\O_K$, using e.g.~\cite[Alg. 6.2.9]{coh}.  Renumber so that $\p_2 = \overline{\p_1}$ and $\p_3 = \overline{\p_4}$.
\item \label{st:gen} Compute generators $\alpha_1$ and $\alpha_2$ for the principal ideals $\p_1\p_3$ and $\p_2 \p_3$, respectively, using e.g.~\cite[Alg.~6.5.10]{coh}.
\item Compute a fundamental unit $u$ of $K_0$ with $\abs{u} > 1$, using e.g.~\cite[Alg. 5.7.1]{coh}.
\item For $i \leftarrow 1,2$, do the following:
	\begin{enumerate}
	\item \label{st:unit} If $\abs{\alpha_i} < \sqrt{p}$, set $\alpha_i \leftarrow \alpha_i u$ until $\abs{\alpha_i} = \sqrt{p}$.  If $\abs{\alpha_i} > \sqrt{p}$, set $\alpha_i \leftarrow \alpha_i u^{-1}$ until $\abs{\alpha_i} = \sqrt{p}$.
	\item Compute the characteristic polynomial $h_i(x)$ of $\alpha_i$, using e.g.~\cite[Prop. 4.3.4]{coh}.  
	\item \label{st:changesign} If $K$ is Galois and $h_1(x) = h_2(-x)$, set $\alpha_2 \leftarrow -\alpha_2$ and $h_2(x) \leftarrow h_2(-x)$.
	\item Set $(n_{i,+1},m_{i,+1}) \leftarrow (p + 1 - \frac{h_i'(0)}{p}, h_i(1))$.  Set $(n_{i,-1},m_{i,-1}) \leftarrow (p + 1 + \frac{h_i'(0)}{p}, h_i(-1))$.
	\end{enumerate}
\item Determine whether the Frobenius endomorphism of $J$ has characteristic polynomial equal to $h_i(\pm x)$ for some $i$:
\begin{enumerate}
	\item \label{st:randompt} Choose a random point $P \in J(\F_p)$ and compute $Q_{j,\tau} = [m_{i,\tau}]P$ for $i \in \{1,2\}$, $\tau \in \{\pm 1\}$.  If none of $Q_{i,\tau}$ is the identity, return {\tt false}.  Otherwise, optionally repeat with another random point $P$.
	\item \label{st:countc} If $J$ passes a certain fixed number of trials of Step \ref{st:randompt},
	compute $\#C(\F_p)$.  If $\#C(\F_p) \neq n_{i,\tau}$ for all $i \in \{1,2\}$, $\tau \in  
	\{\pm 1\}$, return {\tt false}.
	\item \label{st:countj} If $\#C(\F_p) = n_{i,\tau}$, compute $\#J(\F_p)$, using e.g.~Baby Step Giant Step \cite[Alg 5.4.1]{coh}.  If $\#J \neq m_{i,\tau}$ for the same $i,\tau$, return {\tt false}.
\end{enumerate}
\item \label{st:output-pi} If $K$ is Galois, output $S = \{\tau \alpha_1, \tau \overline{\alpha_1}, \tau \alpha_2, \tau \overline{\alpha_2}\}$.  If $K$ is not Galois, output $S = \{\tau \alpha_i, \tau \overline{\alpha_i}\}$, using the $i$ determined in Step \ref{st:countj}.
\item Return {\tt true}.
\end{enumerate}
\end{alg}

\begin{proof}
The proof of \cite[Prop. 4]{el} shows that the ideals $\p_1\p_3$ and $\p_2\p_3$ are principal and the Frobenius endomorphism of $J$ corresponds to a generator of one of these ideals or their complex conjugates.  Furthermore, this generator must have complex absolute value $\sqrt{p}$.  The generators determined in Step \ref{st:gen} are unique up to unit multiple, so Step \ref{st:unit} ensures that the absolute values are $\sqrt{p}$, thus making each $\alpha_i$ unique up to complex conjugation and sign.  

If the Frobenius element corresponds to $\alpha_i$ or $\overline{\alpha_i}$, then $h_i(x)$ is the characteristic polynomial of Frobenius, so we can determine this case by checking whether $\#C(\F_p) = n_{i,+1}$ and $\#J(\F_p) = m_{i,+1}$.  Similarly, if the Frobenius element corresponds to $-\alpha_i$ or $-\overline{\alpha_i}$, then $h_i(-x)$ is the characteristic polynomial of Frobenius, so we can determine this case by checking whether $\#C(\F_p) = n_{i,-1}$ and $\#J(\F_p) = m_{i,-1}$.

If $K$ is Galois 
(with Galois group $C_4$), 
then the ideal $(\alpha_2)$ is equal to $(\alpha_1)^\sigma$ for some $\sigma$ generating the Galois group.  Since complex absolute value squared is the same as the norm from $K$ to its real quadratic subfield $K_0$, $\abs{\alpha_1} = \sqrt{p}$ implies that $\abs{\alpha_1^\sigma} = \sqrt{p}$.  Since $\alpha_1^\sigma$ and $\alpha_2$ both generate $(\alpha_2)$ and have absolute value $\sqrt{p}$, we deduce that $\alpha_1^\sigma = \pm \alpha_2$.  Step \ref{st:changesign} ensures that this sign is positive, so $\alpha_1$ and $\alpha_2$ have the same characteristic polynomial $h_i(x)$, and thus the Frobenius element could be any of the elements output by Step \ref{st:output-pi}.  Since $\Aut(K/\Q)$ is generated by $\sigma$ and $\sigma^2$ is complex conjugation, we have output the $\Aut(K/\Q)$-orbit of the Frobenius element.

If $K$ is not Galois, then the Frobenius element must be either $\alpha_i$ or $\overline{\alpha_i}$.  Since $\Aut(K/\Q)$ in this case consists of only the identity and complex conjugation, Step \ref{st:output-pi} outputs the $\Aut(K/\Q)$-orbit of the Frobenius element.
\end{proof}


\section{Constructing a generating set for $\O_K$}
\label{s:generators}

Given the Jacobian $J$ of a genus $2$ curve over $\F_p$ and a primitive quartic CM field $K$, Algorithm \ref{a:zeta} allows us to determine whether there is some $\pi \in \O_K$ that represents the Frobenius endomorphism of $J$.  Since the complex conjugate $\pibar$ represents the Verschiebung endomorphism, if Algorithm \ref{a:zeta} outputs {\tt true} then we have
\begin{equation}
\label{eq:chain}
\Z[\pi,\pibar] \subseteq \End(J) \subseteq \O_K.
\end{equation}
In this section, we assume we are given a $J/\F_p$ and a $\pi$ such that \eqref{eq:chain} holds, and we wish to determine whether $\End(J) = \O_K$.  

Let $\B$ be a $\Z$-module basis for $\O_K$, and consider the collection of elements $\{ \alpha \in \B \setminus \Z \}$.  Since this collection generates $\O_K$ over $\Z[\pi,\pibar]$, it suffices to determine whether or not each element of the collection is an endomorphism of $J$.  Assuming $K$ satisfies some mild hypotheses, Eisentr\"ager and Lauter give one example of a basis $\B$ that suffices to determine the endomorphism ring \cite[Lemma 6]{el}.  However, the method given in \cite{el} lacks an efficient procedure for testing whether a given $\alpha \in \B$ is an endomorphism of $J$.

In this section, we derive from an arbitrary basis $\B$ a set of generators for $\O_K$ over $\Z[\pi,\pibar]$ that are convenient in the sense that there is an efficient probabilistic algorithm  (Algorithm \ref{a:fields} or Algorithm \ref{a:frob}) for determining whether an element of the set is an endomorphism of $J$.  Our findings are summarized in Proposition \ref{prop:generators}.

We begin by observing that since $K = \Q(\pi)$, any $\alpha \in \O_K$ can be expressed as a polynomial $f \in \Q[\pi]$.  Since $\pi$ satisfies a polynomial of degree $4$ (the characteristic polynomial of Frobenius), $f$ can be taken to have degree $3$.  We may thus write
\begin{equation} 
\label{eq:alpha}
	\alpha = \frac{a_0 + a_1 \pi + a_2 \pi^2 + a_3 \pi^3}{n}
\end{equation}
for some integers $a_0,a_1,a_2,a_3,n$.  We assume that $a_0,a_1,a_2,a_3$ have no common factor with $n$, so that $n$ is the smallest integer such that $n \alpha \in \Z[\pi]$.

\begin{rem}
The LLL lattice reduction algorithm \cite{lll}, as implemented by the MAGMA command {\tt LinearRelation}, finds an expression of the form \eqref{eq:alpha} for any $\alpha \in \O_K$ .  Given as input the sequence $[1,\pi,\pi^2,\pi^3,-\alpha]$, the algorithm outputs a sequence $[a_0,a_1,a_2,a_3,n]$ satisfying the relation \eqref{eq:alpha}.
\end{rem}

The following lemma shows that each $\alpha \in \B \setminus \Z$ can be replaced with a collection of elements that generate the same ring, each with a power of a single prime in the denominator of the expression \eqref{eq:alpha}.

\begin{lemma}
\label{l:euclid}
Let $A \subset B$ be commutative rings with $1$, with $[B:A]$ finite.  Suppose $\alpha \in B$, and let $n$ be the smallest integer such that $n \alpha \in A$.  Suppose $n$ factors into primes as $\ell_1^{d_1} \cdots \ell_r^{d_r}$.  Then
\[ A[\alpha] = \textstyle{A\Big[\frac{n}{\ell_1^{d_1}}\alpha, \ldots, \frac{n}{\ell_r^{d_r}}\alpha\Big].} \]
\end{lemma}

\begin{proof}
Clearly the ring on the right is contained in the ring on the left, so we must show that $\alpha$ is contained in the ring on the right.  It suffices to show that there are integers $c_i$ such that
\begin{equation}
\label{eq:euclid}
 c_1 \frac{n}{\ell_1^{d_1}} + \cdots + c_r \frac{n}{\ell_r^{d_r}} = 1,
\end{equation}
for then we can multiply this identity by $\alpha$ to get our result.  We use the extended Euclidean algorithm and induct on $r$, the number of distinct primes dividing $n$.  If $r = 1$ the result is trivial, for in this case $n/\ell_1^{d_1} = 1$.  Now suppose \eqref{eq:euclid} holds for any $n$ that is divisible by $r$ distinct primes.  If $n'$ is divisible by $r+1$ distinct primes, we can write $n' = n \ell_{r+1}^{d_{r+1}}$ for some $n$ divisible by $r$ distinct primes.  Since $\ell_{r+1}$ is relatively prime to $n$, we can use the extended Euclidean algorithm to write $a\ell_{r+1}^{d_{r+1}} + bn = 1$ for some integers $a,b$.  We can then multiply the first term by the left-hand side of \eqref{eq:euclid} (which is equal to $1$) to get
\begin{equation*}
 a c_1 \frac{n \ell_{r+1}^{d_{r+1}}}{\ell_1^{d_1}} + \cdots + a c_r \frac{n \ell_{r+1}^{d_{r+1}}}{\ell_r^{d_r}} + bn \\
 = ac_1 \frac{n'}{\ell_1^{d_1}} + \cdots + a c_r \frac{n'}{\ell_r^{d_r}}+ b \frac{n'}{\ell_{r+1}^{d_{r+1}}} \\
 = 1.
\end{equation*}
This is an equation of the form \eqref{eq:euclid} for $n'$, which completes the proof.
\end{proof}


The next lemma shows that only primes dividing the index $[\O_K : \Z[\pi]]$ appear in the denominators.
\begin{lemma}
\label{l:dividesindex}
Let $\alpha$ be an element of $\O_K$, and suppose $n$ is the smallest integer such that $n\alpha \in \Z[\pi]$.  Then $n$ divides the index $[\O_K : \Z[\pi]]$.
\end{lemma}

\begin{proof}
Let $N = [\O_K : \Z[\pi]]$.  By definition, $N$ is the size of the abelian group $\O_K/\Z[\pi]$.  Thus we can write any $\alpha \in \O_K$ as $\alpha = a + b$ with $b \in \Z[\pi]$ and $N\cdot a \in \Z[\pi]$.  This shows that $\O_K$ is contained in $\frac{1}{N}\Z[\pi]$.  We may thus write $\alpha = f(\pi)/N$ for a unique polynomial $f$ with integer coefficients and degree at most $3$.  Furthermore, since $n\alpha$ is the smallest multiple of $\alpha$ in $\Z[\pi]$, we may write $\alpha = g(\pi)/n$ for a unique polynomial $g$ with integer coefficients and degree at most $3$, such that $n$ has no factor in common with all the coefficients of $g$.  We thus have $n\cdot f(\pi) = N\cdot g(\pi)$.  If we let $d$ be the gcd of the coefficients of $f$ and $e$ be the gcd of the coefficients of $g$, then we have $n\cdot d = N \cdot e$.  

Let $\ell$ be a prime dividing $e$.  Since $\gcd(n,e) = 1$, $\ell$ must divide $d$, so we can cancel $\ell$ from both sides and get $n \cdot d' = N \cdot e'$ with $e' < e$.  Proceeding in this manner until $e' = 1$, we conclude that $n$ divides $N$.
\end{proof}

We now know that each $\alpha \in \B \setminus \Z$ can be replaced with a collection of elements 
$\{\frac{n}{\ell_i^{d_i}}\alpha\}$, and the only $\ell_i$ appearing are divisors of the index the index $[\O_K : \Z[\pi]]$.  The following lemma and corollary show that
for any $\ell$ which divides $[\O_K : \Z[\pi]]$ exactly (i.e.~$\ell \mid [\O_K: \Z[\pi]]$ and $\ell^2 \nmid [\O_K : \Z[\pi]]$), the element $\frac{n}{\ell}\alpha$ can be replaced by an element of the form $\frac{\pi^k -1 }{\ell}$.  This replacement is useful since by \cite[Fact 10]{el}, determining whether an element of the form $\frac{\pi^k -1}{\ell}$ is an endomorphism is equivalent to testing the field of definition of the $\ell$-torsion. 

\begin{lemma}
\label{l:nofrob}
Let $A \subset B \subset C$ be abelian groups, with $[C:A]$ finite.  Let $\ell$ be a prime, and suppose $\ell$ divides $[C:A]$ and $\ell^2$ does not divide $[C:A]$.  Suppose there is some $\beta \in B$ such that $\beta \not \in A$ and $\ell \beta \in A$.  Then for any $\alpha \in C$ such that $\ell \alpha \in A$, $\alpha \in B$.
\end{lemma}

\begin{proof}
The hypotheses on $[C:A]$ imply that the $\ell$-primary part of $C/A$ (denoted $(C/A)_\ell$) is isomorphic to $\Z/\ell\Z$, so $(B/A)_\ell$ is either trivial or $\Z/\ell\Z$.  The conditions on $\beta$ imply that $\beta$ has order $\ell$ in $B/A$, so $(B/A)_\ell \cong \Z/\ell\Z \cong (C/A)_\ell$, with the isomorphism induced by the inclusion map $B \hookrightarrow C$.  Since $\alpha$ is in the $\ell$-primary part of $C/A$, $\alpha$ must also be in the $\ell$-primary part of $B/A$, so $\alpha \in B$.
\end{proof}

\begin{cor}
\label{c:nofrob}
Suppose $\ell$ divides $[\O_K : \Z[\pi]]$ exactly and $\beta = \frac{\pi^k-1}{\ell} \not \in \Z[\pi]$. Then $\frac{\pi^k-1}{\ell}$ is an endomorphism of $J$ if and only if any $\alpha \in \O_K \setminus \Z[\pi]$ with $\ell \alpha \in \Z[\pi]$ is also an endomorphism.
\end{cor}
\begin{proof}
The result follows directly from Lemma \ref{l:nofrob}, with $A = \Z[\pi]$, $B = \End(J)$, and $C = \O_K$.
\end{proof}

Furthermore, if $p \nmid [\O_K : \Z[\pi,\pibar]]$, then any element $\alpha_i$ with denominator $\ell_i = p$ may be ignored due to the following corollary.  

\begin{cor}
\label{c:pibar}
Suppose $p \nmid [\O_K : \Z[\pi,\pibar]]$.  Then for any $\alpha \in \O_K$ such that $p \alpha \in \Z[\pi]$, $\alpha \in \Z[\pi,\pibar]$.
\end{cor}

\begin{proof}
Since $\pi$ is the Frobenius element, it satisfies a characteristic polynomial of the form
\begin{equation}
\label{eq:frob}
\pi^4 + s_1 \pi^3 + s_2 \pi^2 + s_1 p \pi + p^2 = 0.
\end{equation}
Using $\pi \pibar = p$ and dividing this equation by $\pi$ gives
\begin{equation}
\label{eq:pibar2}
\begin{split}
	\pi^3 + s_1 \pi^2 + s_2 \pi + s_1 p + p \pibar  =  & ~0. \\
\end{split}
\end{equation}
From this equation we see that $p \pibar \in \Z[\pi]$, so either $[\Z[\pi,\pibar] : \Z[\pi]] = p$ or $\pibar \in \Z[\pi]$.  If $\pibar \in \Z[\pi]$ then $p$ divides the coefficients of the terms on the left hand side of \eqref{eq:pibar2}, which it does not, so we deduce that $\pibar \not \in \Z[\pi]$ and $[\Z[\pi,\pibar] : \Z[\pi]] = p$. 
The hypothesis $p \nmid [\O_K : \Z[\pi,\pibar]]$ thus implies that $p$ divides $[\O_K : \Z[\pi]]$ exactly, so we may apply Lemma \ref{l:nofrob} with $\ell = p$, $A = \Z[\pi]$, $B = \Z[\pi,\pibar]$, $C = \O_K$, and $\beta = \pibar$.  
\end{proof}

Thus any $\alpha$ satisfying the conditions of the corollary is automatically an endomorphism.  We now show that the condition $p \nmid [\O_K : \Z[\pi,\pibar]]$ is automatically satisfied for all primes $p$ except possibly $2$ and $3$.  

\begin{prop}
\label{p:no-p}
Suppose $p > 3$ and that $\pi \in \O_K$ corresponds to the Frobenius endomorphism of an ordinary abelian surface $A$ over $\F_p$.  Then $p \nmid [\O_K : \Z[\pi,\pibar]]$.
\end{prop}

\begin{proof}  
Let $\Delta(R)$ denote the discriminant of a $\Z$-module $R$.
Christophe Ritzenthaler pointed out that this proposition follows from \cite[Proposition 9.4]{Howe}, which
shows that $$\Disc(\Z[\pi,\pibar]) = \pm \Norm_{K/\Q}(\pi-\pibar) \Disc(\Z[\pi+\pibar]).$$
Alternatively, it is shown in \cite[Proposition 7.4]{LPP} that any prime that divides the index $[\O_K : \Z[\pi,\pibar]]$ must divide
either $[\O_{K_0}:\Z[\pi+\pibar]]$ or $\frac{\Disc(\O_{K_0}[\pi])}{\Disc(\O_K)}$, and, using \cite[Theorem 1.3]{Howe}, that the second quantity is prime to $p$ if the abelian surface is ordinary.  The same proposition also shows that $\Disc(\Z[\pi+\pibar]) < 16p$, and since 
$$\frac{\Disc(\Z[\pi+\pibar])}{\Disc(\O_{K_0})} = [\O_{K_0}:\Z[\pi+\pibar]]^2,$$
we conclude that if $p$ divides $[\O_K : \Z[\pi,\pibar]]$ 
then $p^2$ divides $[\O_{K_0}:\Z[\pi+\pibar]]^2$, and thus
$$p^2 \le \frac{\Disc(\Z[\pi+\pibar])}{\Disc(\O_{K_0})} < \frac{16p}{5}$$ (since a real quadratic field has discriminant at least $5$), which implies $p \le 3$.\end{proof}

The following proposition summarizes the results of this section. 

\begin{prop} \label{prop:generators}
Suppose $\{\alpha_i\}$ generates $\O_K$ as a $\Z$-algebra.  Let $n_i$ be the smallest integer such that $n_i \alpha_i \in \Z[\pi]$, and write the prime factorization of $n_i$ as $n_i = \prod_j \ell_{ij}^{d_{ij}}$.  For each $(i,j)$ with $\ell_{ij} \neq p$, let $k_{ij}$ be an integer such that $\pi^{k_{ij}} - 1 \in \ell_{ij} \O_K$.  Suppose $p > 3$.  Then the following set generates $\O_K$ over $\Z[\pi,\pibar]$:
\begin{equation*}
\left\{ \frac{n_i}{\ell_{ij}^{d_{ij}}} \alpha_i : 
	\ell_{ij}^2 \mid [\O_K : \Z[\pi]] \right\} 
\cup \left\{ \frac{\pi^{k_{ij}}-1}{\ell_{ij}} : 
	\ell_{ij}^2 \nmid [\O_K : \Z[\pi]], \ell_{ij} \neq p \right\}.
\end{equation*}
\end{prop}

\begin{rem} Proposition \ref{prop:generators} shows that if $p > 3$ and the index $[\O_K : \Z[\pi,\pibar]]$ is square-free, then $\O_K$ can be generated over $\Z[\pi,\pibar]$ by a collection of elements of the form $\frac{\pi^k-1}{\ell}$.  This answers a question raised by Eisentr\"ager and Lauter \cite[Remark 5]{el}.
\end{rem}

In our application, $\pi \in \O_K$ is only determined up to an automorphism of $K$, but Proposition \ref{prop:generators} can still be used to determine a generating set for $\O_K$.

\begin{cor}
\label{c:galois-generators}
Let $\mathcal{S} \subset \O_K$ be the set given in Proposition \ref{prop:generators}.  Let $\sigma$ be an element of $\Aut(K/\Q)$.  Then the set $\{\beta^\sigma : \beta \in \mathcal{S}\}$ generates $\O_K$ over $\Z[\pi^\sigma, \pibar^\sigma]$.
\end{cor}

\begin{proof}
By Proposition \ref{prop:generators}, the set $\{\pi,\pibar\} \cup \mathcal{S}$ generates $\O_K$ as a $\Z$-algebra.  Since $\O_K$ is mapped to itself by $\Aut(K/\Q)$, the set $\{\pi^\sigma,\pibar^\sigma\} \cup \{\beta^\sigma : \beta \in \mathcal{S}\}$ also generates $\O_K$ as a $\Z$-algebra.  The statement follows immediately.
\end{proof}


\section{Determining fields of definition}
\label{s:fields}

In this section, we consider the problem of determining the field of definition of the $n$-torsion points of the Jacobian $J$ of a genus 2 curve over $\F_p$.  By \cite[Fact 10]{el}, the $n$-torsion points of $J$ are defined over $\F_{p^k}$ if and only if $(\pi^k-1)/n$ is an endomorphism of $J$, where $\pi$ is the Frobenius endomorphism of $J$.  Thus determining the field of definition of the $\ell$-torsion points allows us to determine whether some of the elements given by Proposition \ref{prop:generators} are endomorphisms.

\begin{alg}
\label{a:torsdef}
The following algorithm takes as input a primitive quartic CM field $K$, an element $\pi \in \O_K$ with $\pi\pibar = p$, and an integer $n$ with $\gcd(n,p) = 1$, and outputs the smallest integer $k$ such that $\pi^k-1 \in n\O_K$.  If $J$ is the Jacobian of a genus 2 curve over $\F_p$ with Frobenius $\pi^\sigma$ for some $\sigma \in \Aut(K/\Q)$ and $\End(J) = \O_K$, this integer $k$ is such that the $n$-torsion points of $J$ are defined over $\F_{p^k}$.

\begin{enumerate}
\item Compute a $\Z$-basis $\B = (1,\delta,\gamma,\kappa)$ of $\O_K$, using \cite{SpeWi} or \cite[Algorithm 6.1.8]{coh}, and write $\pi = (a,b,c,d)$ in this basis.  Set $k \leftarrow 1$.
\item Let $\bar \B$ be the reduction of the elements of $\B$ modulo $n$.  Let $(a_1,b_1,c_1,d_1) = (a,b,c,d) \pmod{n}$.
\item \label{st:pik} Compute $\pi^k \equiv (a_k,b_k,c_k,d_k) \pmod{n}$ with respect to $\bar \B$.
\item \label{st:gcd} If $(a_k,b_k,c_k,d_k) \equiv (1,0,0,0) \pmod{n}$, output $k$.  Otherwise set $k \leftarrow k+1$ and go to Step \ref{st:pik}.
\end{enumerate}
\end{alg}

\begin{proof}
The set $\bar \B$ is a $\Z/n\Z$-basis of $\O_K/n\O_K$, so if $\pi^k \equiv (1,0,0,0) \pmod n$, then $\pi^k-1 \in n\O_K$ (since the first element of $\bar \B$ is $1$).  
Since $n\O_K$ is mapped to itself by $\Aut(K/\Q)$, we have $(\pi^\sigma)^k-1 \in n\O_K$.  If $\End(J) = \O_K$, then $\frac{(\pi^\sigma)^k-1}{n} \in \O_K = \End(J)$, so by \cite[Fact 10]{el}, $J[n] \subset J(\F_{p^k})$.
\end{proof}

\begin{rem}
Since $J[n] = \bigoplus J[\ell^d]$ for prime powers $\ell^d$ dividing $n$, we may speed up Algorithm \ref{a:torsdef} by factoring $n$ and computing $k(\ell^d)$ for each prime power factor $\ell^d$; then $k(n) = \lcm(k(\ell^d))$.  Furthermore, we will see in Propositions \ref{p:ltors} and \ref{p:ldtors} below that for a fixed $\ell^d$, the possible values of $k$ are very limited.  Thus we may speed up the algorithm even further by precomputing these possible values and testing each one, rather than increasing the value of $k$ by $1$ until the correct value is found.
\end{rem}

Eisentr\"ager and Lauter \cite{el} computed endomorphism rings in several examples by determining the group structure of $J(\F_{p^k})$ to decide whether $J[n] \subset J(\F_{p^k})$.  This is an exponential-time algorithm that is efficient only for very small $k$.  Eisentr\"ager and Lauter also suggested that the algorithm of Gaudry-Harley \cite{gh} could be used to determine the field of definition of the $n$-torsion points.  
One of the primary purposes of this article is to present an efficient probabilistic algorithm to test the field of definition of $J[n]$.  
Below we describe the various methods of testing the field of definition of the $n$-torsion of $J$.
Since $J[n] = \bigoplus J[\ell^d]$ as $\ell^d$ ranges over maximal prime-power divisors of $n$, it suffices to consider each prime-power factor separately.  We thus assume in what follows that $n = \ell^d$ is a prime power.

\subsection{The brute force method}
\label{ss:brute1}

The simplest method of determining the field of definition of the $n$-torsion is to compute the abelian group structure of $J(\F_{p^k})$.  The MAGMA syntax for this computation is straightforward, and the program returns a group structure of the form
$$ J(\F_{p^k}) \cong \frac{\Z}{a_1\Z} \times \cdots \times \frac{\Z}{a_j\Z}, $$
with $a_1 \mid \cdots \mid a_j$.  The $n$-torsion of $J$ is contained in $J(\F_{p^k})$ if and only if $j = 4$ and $n$ divides $a_1$.

While this method is easy to implement, if $k$ is too large it may take too long to compute the group structure (via Baby-Step/Giant-Step or similar algorithms), or even worse we may not even be able to factor $\#\Jac(C)(\F_{p^k})$.  In practice, computing group structure in MAGMA seems to be feasible for group sizes up to roughly $2^{200}$, which means $p^k$ should be no more than roughly $2^{100}$, and thus $k$ will have to be very small.  Thus the brute force method is very limited in scope; however, it has the advantage that in the small cases it can handle it runs fairly quickly and always outputs the right answer.

\subsection{The Gaudry-Harley-Schost method}

\label{ss:gaudry}

Gaudry and Harley \cite{gh} define a Schoof-Pila-like algorithm for counting points on genus $2$ curves.  The curves input to this algorithm are assumed to have a degree $5$ model over $\F_q$, so we can write elements of the Jacobian as pairs of affine points minus twice the Weierstrass point at infinity.  An intermediate step in the algorithm is to construct a polynomial $R(x) \in \F_q[x]$ with the following property: if $P_1$ and $P_2$ are points on $C$ such that $D = [P_1] + [P_2] - 2 [\infty]$ is an $n$-torsion point of $J$, then the $x$-coordinates of $P_1$ and $P_2$ are roots of $R$.  
The field of definition of the $x$-coordinates is at most a degree-two extension of the field of definition of $D$.  
Thus in many cases the field of definition of the $n$-torsion points can be determined from the factorization of $R(x)$.  

Gaudry has implemented the algorithm in MAGMA and NTL; the algorithm involves taking two resultants of pairs of two-variable polynomials of degree roughly $n^2$.  The algorithm uses the clever trick of computing a two-variable resultant by computing many single-variable resultants and interpolating the result.  The interpolation only works if the field of definition of $J$ has at least $4n^2 - 8n + 4$ elements, so we must base extend $J$ until the field of definition is large enough.  Since $R(x)$ has coefficients in $\F_p$, this base extension has no effect on the result of the computation.

Gaudry and Harley's analysis of the algorithm gives a running time of $\tilde{O}(n^6)$ field multiplications if fast polynomial arithmetic is used, and $O(n^8)$ otherwise.  
Due to its large space requirements, the algorithm has
only succeeded at handling inputs of size $n \le 19$ \cite{GS}.

\subsection{A probabilistic method}
\label{ss:prob1}

As usual, we let $J$ be the Jacobian of a genus 2 curve over $\F_{p^k}$, and $\ell \neq p$ be a prime.  Let $H$ be the $\ell$-primary part of $J(\F_{p^k})$.  Then $H$ has the structure
$$ H = \frac{\Z}{\ell^{\alpha_1}\Z} \times \frac{\Z}{\ell^{\alpha_2}\Z} \times \frac{\Z}{\ell^{\alpha_3}\Z} \times \frac{\Z}{\ell^{\alpha_4}\Z},$$
with $\alpha_1 \leq \alpha_2 \leq \alpha_3 \leq \alpha_4$.  Our test rests on the following observations:
\begin{itemize}
\item If the $\ell^d$-torsion points of $J$ are defined over $\F_{p^k}$, then $\alpha_1 \geq d$, and the number of $\ell^d$-torsion points in $H$ is $\ell^{4d}$.
\item If the $\ell^d$-torsion points of $J$ are not defined over $\F_{p^k}$, then $\alpha_1 < d$, and the number of $\ell^d$-torsion points in $H$ is at most $\ell^{4d-1}$.
\end{itemize}We thus make the following calculation: write $\#J(\F_{p^k}) = \ell^s m$ with $\ell \nmid m$.  Choose a random point $P \in J$.  Then $[m]P \in H$, and we test whether $[\ell^d m]P = O$ in $J$.  If the $\ell^d$-torsion points of $J$ are defined over $\F_{p^k}$, then $[\ell^d m]P = O$ with probability $\rho = \ell^{4d-s}$, while if the $\ell^d$-torsion points of $J$ are not defined over $\F_{p^k}$ then $[\ell^d m]P = O$ with probability at most $\rho/\ell$.  If we perform the test enough times, we can determine which probability distribution we are observing and thus conclude, with a high degree of certainty, whether the $\ell^d$-torsion points are defined over $\F_{p^k}$.

This method is very effective in practice, and can be implemented for large $k$: while computing the group structure of $J(\F_{p^k})$ for large $k$ may be infeasible, it is much easier to compute {\it points} on $J(\F_{p^k})$ and to do arithmetic on those points.  We now give a formal description of the algorithm and determine its probability of success.

\begin{alg}
\label{a:fields}
The following algorithm takes as input the Jacobian $J$ of a genus 2 curve defined over a finite field $\F_{q}$, a prime power $\ell^d$ with $\gcd(\ell,q) = 1$, and a real number $\epsilon \in (0,1)$.  If $J[\ell^d] \subset J(\F_q)$, then the algorithm outputs {\tt true} with probability at least $1 - \epsilon$.  If $J[\ell^d] \not \subset J(\F_q)$, then the algorithm outputs {\tt false} with probability at least $1 - \epsilon$.
\begin{enumerate}
\item Compute $\#J(\F_q) = \ell^s m$, where $\ell \nmid m$.  If $s < 4d$ output {\tt false}.
\item Set $\rho \leftarrow \ell^{4d-s}$, $N \leftarrow \ceil{\frac{\sqrt{-2 \log\epsilon}}{\rho} ( \frac{2\ell}{\ell-1})}$, $B \leftarrow \rho N \left( \frac{\ell+1}{2\ell} \right)$.
\item \label{st:fieldloop} Repeat $N$ times:
	\begin{enumerate}
	\item Choose a random point $P_i \in J(\F_q)$.
	\item Compute $Q_i \leftarrow [\ell^d m]P_i$
	\end{enumerate}
\item If at least $B$ of the $Q_i$ are the identity element $O$ of $J$, output {\tt true}; otherwise output {\tt false}.
\end{enumerate}
\end{alg}

\begin{proof}
As observed above, if $J[\ell^d] \subset J(\F_q)$, then $Q_i = O$ with probability $\rho$, while if $J[\ell^d] \not \subset J(\F_q)$, then $Q_i = O$ with probability at most $\rho/\ell$.  Thus all we have to do is compute enough $Q_i$ to distinguish the two probability distributions.  To figure out how many ``enough'' is, we use the Chernoff bound \cite[Ch.~8, Prop.~5.3]{ross}.  The version of the bound we use is as follows: If $N$ weighted coins are flipped and $\mu$ is the expected number of heads, then for any $\delta \in (0,1]$ we have
\begin{equation}
\label{eq:chernoff}
\begin{split}
\Pr [\mbox{\#heads} < (1-\delta )\mu ] & < e^{-\mu^2 \delta^2/2} \\
\Pr [\mbox{\#heads} > (1+\delta )\mu ] & < e^{-\mu^2 \delta^2/2}.
\end{split}
\end{equation}
In our case we are given two different probability distributions for the coin flip and wish to tell them apart.  If the $\ell^d$-torsion points of $J$ are defined over $\F_q$, then the probability that $Q_i = O$ is $\rho = \ell^{4d}/\ell^s$.  Thus the expected number of $Q_i$ equal to $O$ is $\mu_1 = \rho N$.  If the $\ell^d$-torsion points are not defined over $\F_q$, then the expected number of $Q_i$ equal to $O$ is at most $\mu_2 = \rho N/\ell$.  Thus if we set $B = \rho N (\frac{\ell+1}{2\ell})$ to be the midpoint of $[\mu_2,\mu_1]$, we will deduce that $J[\ell^d] \subset J(\F_q)$ if the number of $Q_i$ equal to $O$ is at least $B$, and $J[\ell^d] \not \subset J(\F_q)$ otherwise.

We thus wish to find an $N$ such that this deduction is correct with probability at least $1-\epsilon$, i.e.~an $N$ such that
\begin{equation}
\label{eq:prob}
\begin{split}
\Pr[\#\{Q_i : Q_i = O \} < B] & <  \epsilon \quad \mbox{if $J[\ell^d] \subset J(\F_q)$}, \\
\Pr[\#\{Q_i : Q_i = O \} > B] & <  \epsilon \quad \mbox{if $J[\ell^d] \not \subset J(\F_q)$}.
\end{split}
\end{equation}
Substituting our choice of $B$ into the Chernoff bound \eqref{eq:chernoff} gives
\begin{equation}
\begin{split}
\Pr[\#\{Q_i : Q_i = O \} < B] & < e^{-2\mu_1^2 (\frac{\ell-1}{4\ell})^2} \quad \mbox{if $J[\ell^d] \subset J(\F_q)$}, \\
\Pr[\#\{Q_i : Q_i = O \} > B] & < e^{-2\mu_2^2 (\frac{\ell-1}{4})^2} \quad \mbox{if $J[\ell^d] \not \subset J(\F_q)$}.
\end{split}
\end{equation}
From these equations, we see that we wish to have $2\mu_1^2 (\frac{\ell-1}{4\ell})^2 > -\log \epsilon$ and $2\mu_2^2 (\frac{\ell-1}{4})^2 > - \log \epsilon$.  The two left sides are equal since $\mu_2 = \mu_1/\ell$.  We thus substitute $\mu_1 = \rho N$ into the relation $2\mu_1^2 (\frac{\ell-1}{4\ell})^2 > -\log \epsilon$, and find that
\begin{equation*}
N > \frac{\sqrt{-2 \log\epsilon}}{\rho} \cdot \frac{2\ell}{\ell-1}.
\end{equation*}
Thus this value of $N$ suffices to give the desired success probabilities.
\end{proof}

\begin{rem}
If $s = 4d$, then the algorithm can be simplified considerably.  In this case, if $J[\ell^d] \subset J(\F_q)$ then the $\ell$-primary part $H$ of $J(\F_q)$ is isomorphic to $(\Z/\ell^d\Z)^4$, and if not then it contains a point of order greater than $\ell^d$.  Thus if $J[\ell^d] \subset J(F)$ then $Q_i$ will always be the identity, and the algorithm will always return {\tt true}.  On the other hand, if $J[\ell^d] \not \subset J(\F_q)$, we may abort the algorithm and return {\tt false} as soon as we find a point $Q_i \ne O$, for in this case we have found a point in $H$ of too large order, and thus the $\ell^d$-torsion points are not defined over $\F_q$.  If $J[\ell^d] \not \subset J(\F_q)$, then the probability that a random point in $H$ has order $\leq \ell^d$ is at most $1/\ell$, so we must conduct at least $N = \ceil{\frac{-\log \epsilon}{\log \ell}}$ trials to ensure a success probability of at least $1- \epsilon$. Thus in this case the method may require many fewer trials.
\end{rem}

\begin{rem}
Note that while $\#J(\F_q)$ may be very large, in our application where $J$ is defined over a small prime field it is easy to compute $\#J(\F_q)$ from the zeta function of the curve of which $J$ is the Jacobian.  Furthermore, while it is probably impossible to factor $\#J(\F_q)$ completely in a reasonable amount of time, it is easy to determine the highest power of $\ell$ that divides $\#J(\F_q)$.
\end{rem}

\begin{prop}
\label{p:fieldtime}
Let $J$ be the Jacobian of a genus 2 curve over $\F_p$.  
Assume that the zeta function of $J/\F_p$ is known, so that the cost to compute $\#J(\F_{p^k}) = \ell^s m$ is negligible.
Then the expected number of operations in $\F_p$ necessary to execute Algorithm \ref{a:fields} on $J/\F_{p^k}$ (ignoring $\log \log p$ factors) is 
\begin{equation*}
O(k^2  \log k (\log^2 p) \ell^{s - 4d}(-\log \epsilon)^{1/2})
\end{equation*}
\end{prop}

\begin{proof}
We must compare the cost of the two actions of Step \ref{st:fieldloop}, repeated $N$ times.  Choosing a random point on $J(\F_q)$ is equivalent to computing a constant number of square roots in $\F_q$, and taking a square root requires $O(\log q)$ field operations in $\F_q$ (see \cite[Algorithm 14.15 and Corollary 14.16]{gg}).  The order of $J(\F_q)$ is roughly $q^2$, so multiplying a point on $J(\F_q)$ by an integer using a binary expansion takes $O(\log q)$ point additions on $J(\F_q)$.  Each point addition takes a constant number of field operations in $\F_q$, so we see that the time of each trial is $O(\log q) = O(k \log p)$.  If fast multiplication techniques are used, then the number of field operations in $\F_p$ needed to perform one field operation in $\F_q$ is $O(\log q \log \log q) = O(k \log k \log p)$ (ignoring $\log \log p$ factors), so each trial takes $O(k^2 \log k \log^2 p)$ field operations in $\F_p$.  The number of trials is $O(\ell^{s-4d} \sqrt{-\log \epsilon})$, which gives a total of $O(k^2  \log k (\log^2 p) \ell^{s - 4d}(-\log \epsilon)^{1/2})$ field operations in $\F_p$.
\end{proof}



\section{Computing the action of Frobenius}
\label{s:frob}

As in the previous section, we consider a genus 2 curve $C$ over $\F_p$ with Jacobian $J$, and assume that the endomorphism ring of $J$ is an order in the ring of integers $\O_K$ of a primitive quartic CM field $K$.  We let $\pi$ represent the Frobenius endomorphism, and we look at elements $\alpha \in \O_K$ such that $\ell^d \alpha \in \Z[\pi]$ for some prime power $\ell^d$.  We wish to devise a test that, given such an $\alpha$, determines whether $\alpha$ is an endomorphism of $J$.

Since $\pi$ satisfies a quartic polynomial with integer coefficients, we can write $\alpha$ as
\begin{equation}
\label{eq:alpha-l}
\alpha = \frac{a_0 + a_1 \pi + a_2 \pi^2 + a_3 \pi^3}{\ell^d}
\end{equation}
for some integers $a_0,a_1,a_2,a_3$.
Expressing $\alpha$ in this form is useful because of the following fact proved by Eisentr\"ager and Lauter \cite[Corollary 9]{el}: $\alpha$ is an endomorphism if and only if $T = a_0 + a_1 \pi + a_2 \pi^2 + a_3 \pi^3$ acts as zero on the $\ell^d$-torsion.  Thus we need a method for determining whether $T$ acts as zero on the $\ell^d$-torsion.  Since $T$ is a linear operator, it suffices to check whether $T(Q_i)$ is zero for each $Q_i$ in some set whose points span the full $\ell^d$-torsion.  Below we describe three different ways to compute such a spanning set.

\subsection{The brute force method}
\label{ss:brute2}

The most straightforward way to compute a spanning set for the $\ell^d$-torsion is to use group structure algorithms to compute a basis of $J[\ell^d]$.  This method was used in \cite{el} to compute the class polynomials in one example.  The methods of Section \ref{s:fields} determine a $k$ for which $J[\ell^d] \subset J(\F_{p^k})$.  The computation of the group structure of $J(\F_{p^k})$ gives generators for the group; multiplying these generators by appropriate integers gives generators for the $\ell^d$-torsion.  It is then straightforward to compute the action of $T$ on each generator $g_i$ for $1 \leq i \le 4$.  If $T(g_i) = O$ for all $i$, then $\alpha$ is an endomorphism; otherwise $\alpha$ is not an endomorphism.

This method of computing a spanning set has the same drawback as the brute-force method of computing fields of definition: since the best algorithm for computing group structure runs in time exponential in $k \log p$, the method becomes prohibitively slow as $k$ increases.  Thus the method is only effective when $\ell^d$ is very small.  

\subsection{A probabilistic method}
\label{ss:prob2}

The method of Section \ref{ss:brute2} for computing generators of $J[\ell^d]$ becomes prohibitively slow as the field of definition of the $\ell^d$-torsion points becomes large.  However, we can get around this obstacle by randomly choosing many points $Q_i$ of exact order $\ell^d$, so that it is highly probable that the set $\{Q_i\}$ spans $J[\ell^d]$.

Recall that we wish to test whether the operator $T = a_0 + a_1 \pi + a_2 \pi^2 + a_3 \pi^3$ acts as zero on the $\ell^d$-torsion.  To perform the test, we determine the field $\F_{p^k}$ over which we expect the $\ell^d$-torsion to be defined.  (See Section \ref{s:fields}.)  We pick a random point $P \in J(\F_{p^k})$ and multiply $P$ by an appropriate integer to get a point $Q$ whose order is a power of $\ell$.  If $Q$ has order $\ell^d$, we act on $Q$ by the operator $T$ and test whether we get the identity of $J$; otherwise we try again with a new $P$.  (See Section \ref{ss:couveignes} for another method of randomly choosing $\ell^d$-torsion points.)  We repeat the test until it is overwhelmingly likely that the points $Q$ span the $\ell^d$-torsion.  If the set of $Q$ spans the $\ell^d$-torsion, then $\alpha$ is an endomorphism if and only if $T$ acts as zero on all the $Q$.

\begin{alg}
\label{a:frob}
The following algorithm takes as input the Jacobian $J$ of a genus 2 curve over $\F_{q}$ with CM by $K$, a prime power $\ell^d$ with $\gcd(\ell,q) = 1$, the element $\pi \in \O_K$ corresponding to the Frobenius endomorphism of $J$, an element $\alpha \in \O_K$ such that $\ell^d \alpha \in \Z[\pi]$, and a real number $\epsilon > 0$.  The algorithm outputs {\tt true} or {\tt false}.

Suppose $J[\ell^d] \subset J(\F_{q})$.  If $\alpha$ is an endomorphism of $J$, then the algorithm outputs {\tt true}.  If $\alpha$ is not an endomorphism of $J$, then the algorithm outputs {\tt false} with probability at least $1 - \epsilon$.

\begin{enumerate}
\item \label{st:rep} Compute $a_0,a_1,a_2,a_3$ such that $\alpha$ satisfies equation \eqref{eq:alpha-l}.
\item \label{st:frbound} Set $N$ to be
	\begin{equation*}
	N \leftarrow \left\{ \begin{array}{cc}
	\ceil{\frac{1}{d - \log_\ell 2}\left(-\log_\ell \epsilon + 3d\right)} & \mbox{if $\ell^d > 2$} \\
	\max \{ \ceil{-2 \log_2 \epsilon} + 6, 16\} & \mbox{if $\ell^d = 2$}. 
	\end{array} \right.
	\end{equation*}
\item \label{st:count} Compute $\#J(\F_{q}) = \ell^s m$, where $\ell \nmid m$. 
\item Set $i \leftarrow 1$.  
\item \label{st:rand} Choose a random point $P_i \in J(\F_{q})$.  Set $Q_i \leftarrow [m]P_i$.  Repeat until $[\ell^d] Q_i = O$ and $[\ell^{d-1}] Q_i \neq O$.
\item \label{st:check} Compute 
	\begin{equation}
	\label{eq:check}
	[a_0]Q_i + [a_1]\Frob_p(Q_i) + [a_2]\Frob_{p^2}(Q_i) + [a_3]\Frob_{p^3}(Q_i)
	\end{equation}
in $J(\F_{q})$.  If the result is nonzero output {\tt false}.
\item If $i < N$, set $i \leftarrow i+1$ and go to Step \ref{st:rand}.
\item Output {\tt true}.
\end{enumerate}
\end{alg}

\begin{proof}
By \cite[Corollary 9]{el}, $\alpha$ is an endomorphism of $J$ if and only if the expression \eqref{eq:check} is $O$ for all $\ell^d$-torsion points $Q$.  Furthermore, it suffices to check the the expression only on a basis of the $\ell^d$-torsion.  Step \ref{st:rand} repeats until we find a point $Q_i$ of exact order $\ell^d$; the assumption $J[\ell^d] \subset J(\F_q)$ guarantees that we can find such a point.  The algorithm computes a total of $N$ such points $Q_i$.  Thus if the set of $Q_i$ span $J[\ell^d]$, then the algorithm will output {\tt true} or {\tt false} correctly, according to whether $\alpha \in \End(J)$.  We must therefore compute a lower bound for the probability that the set of $Q_i$ computed span $J[\ell^d]$.

To compute this bound, we will compute an upper bound for the probability that $N$ points of exact order $\ell^d$ do not span $J[\ell^d]$.  We will make repeated use of the following inequality, which can be proved easily with simple algebra: if $\ell$, $d$, $n$, and $m$ are positive integers with $\ell > 1$ and $n > m$, then
\begin{equation}
\frac{\ell^{md} - \ell^{m(d-1)}}{\ell^{nd} - \ell^{n(d-1)}} < \frac{1}{\ell^{(n-m)d}} 
\end{equation}

Next we observe that in any group of the form $(\Z/\ell^d\Z)^r$, there are $\ell^{rd} - \ell^{r(d-1)}$ elements of exact order $\ell^d$.  The probability that a set of $N$ elements does not span a $4$-dimensional space is the sum of the probabilities that all the elements span a $j$-dimensional subspace, for $j = 1,2,3$.  We consider each case:

\begin{itemize}
\item $j = 1$: All of the $Q_i$ are in the space spanned by $Q_1$, and $Q_1$ can be any element.  The probability of this happening is
$$ \left(\frac{\ell^d - \ell^{d-1}}{\ell^{4d} - \ell^{4(d-1)}}\right)^{N-1} 
< \left(\frac{1}{\ell^{3d}}\right)^{N-1}. $$
\item $j = 2$: $Q_1$ can be any element, one of the $Q_i$ must be independent of $Q_1$, and the remaining $N-2$ elements must be in the same $2$-dimensional subspace.  There are $N-1$ ways to choose the second element, so the total probability is
$$ (N-1)\left(1 - \frac{\ell^d - \ell^{d-1}}{\ell^{4d} - \ell^{4(d-1)}}\right)\left(\frac{\ell^{2d} - \ell^{2(d-1)}}{\ell^{4d} - \ell^{4(d-1)}}\right)^{N-2} < N \left(\frac{1}{\ell^{2d}}\right)^{N-2}. $$
\item $j = 3$: $Q_1$ can be any element, and there must be two more linearly independent elements; there are $\binom{N-1}{2}$ ways of choosing these elements.  The remaining $N-3$ elements must all be in the same $3$-dimensional subspace, so the total probability is
\begin{equation*}
\begin{split}
\frac{(N-1)(N-2)}{2}\left(1 - \frac{\ell^d - \ell^{d-1}}{\ell^{4d} - \ell^{4(d-1)}}\right)
\left(1 - \frac{\ell^{2d} - \ell^{2(d-1)}}{\ell^{4d} - \ell^{4(d-1)}}\right)
\left(\frac{\ell^{3d} - \ell^{3(d-1)}}{\ell^{4d} - \ell^{4(d-1)}}\right)^{N-3} \\
< \frac{N^2}{2} \left(\frac{1}{\ell^{d}}\right)^{N-3}. 
\end{split}
\end{equation*}
\end{itemize}
Summing these three cases, we see that the total probability that the $Q_i$ do not span $J[\ell^d]$ is bounded above by
\begin{equation}
\label{eq:bound1}
 N^2 \left(\frac{1}{\ell^d}\right)^{N-3}. 
\end{equation}
Since $2^N \geq N^2$ for $N \geq 4$, we have
\begin{equation*}
 N^2 \left(\frac{1}{\ell^d}\right)^{N-3} \le \ell^{-dN + 3d + N \log_{\ell} 2}.
\end{equation*}
(Note that $N \geq 4$ must always hold if we want to have a spanning set of $J[\ell]$.)  Setting this last expression less than $\epsilon$ and taking logs, we find 
\begin{equation}
\label{eq:frbound}
N \ge \frac{1}{d - \log_\ell 2}\left(-\log_\ell \epsilon + 3d\right).
\end{equation}
Thus if the number of trials $N$ is greater than or equal to the right hand side of \eqref{eq:frbound}, then the probability of success is at least $1 - \epsilon$.  

The right hand side of expression \eqref{eq:frbound} is undefined if $\ell = 2$, $d = 1$, so we must make a different estimate.  Since $2^{N/2} \geq N^2$ for $N \geq 16$, the estimate \eqref{eq:bound1} bounds the probability of $Q_i$ not spanning $J[\ell^d]$ by
\begin{equation*}
\frac{N^2}{2^{N-3}} \le \frac{1}{2^{N/2-3}}.
\end{equation*}
Setting the right hand side less than $\epsilon$ and taking logs gives
\begin{equation}
\label{eq:fr2bound}
N \ge -2 \log_2 \epsilon + 6.
\end{equation}
Thus if the number of trials $N$ is greater than or equal to the maximum of $16$ and the right hand side of \eqref{eq:fr2bound}, then the probability of success is at least $1 - \epsilon$.
\end{proof}

\begin{cor}
\label{c:galois-frob}
Let $J$, $\ell^d$, $\alpha$, and $\epsilon$ be as in Algorithm \ref{a:frob}.  Suppose $\pi \in \O_K$ is such that $\pi^\sigma$ corresponds to the Frobenius endomorphism of $J$ for some $\sigma \in \Aut(K/\Q)$.  Suppose $J[\ell^d] \subset J(\F_{q})$, and suppose Algorithm \ref{a:frob} is run with inputs $J$, $\F_q$, $\pi$, $\alpha$, $\epsilon$.  If $\alpha^\sigma$ is an endomorphism of $J$, then the algorithm outputs {\tt true}.  If $\alpha^\sigma$ is not an endomorphism of $J$, then the algorithm outputs {\tt false} with probability at least $1 - \epsilon$.
\end{cor}

\begin{proof}
If we write $\alpha$ in the form \eqref{eq:alpha-l}, then we have \begin{equation} \label{eq:alpha-sigma} \alpha^\sigma = \frac{a_0 + a_1 \pi^\sigma + a_2 (\pi^\sigma)^2 + a_3 (\pi^\sigma)^3}{\ell^d}.
\end{equation}
Step \ref{st:check} of the algorithm determines whether the numerator of this expression acts as zero on $\ell^d$-torsion points.  By \cite[Corollary 9]{el}, this action is identically zero if and only if $\alpha^\sigma$ is an endomorphism of $J$.  The statement now follows from the correctness of Algorithm \ref{a:frob}.
\end{proof}

\begin{rem}
Since $Q_i$ is an $\ell^d$-torsion point in Step \ref{st:check}, we may speed up the computation of the expression \eqref{eq:check} by replacing each $a_j$ with a small representative of $a_j$ modulo $\ell^d$.  We may also rewrite the expression \eqref{eq:check} as
$$ [a_0]Q_i + \Frob_p([a_1]Q_i + \Frob_p([a_2]Q_i + \Frob_p([a_3]Q_i))) $$
to reduce the number of $\Frob_p$ operations from $6$ to $3$.
\end{rem}

\begin{rem} Algorithm \ref{a:frob} assumes that the $\ell^d$-torsion points of $J$ are defined over $\F_{q}$, so with enough trials we are almost certain to get a spanning set of points $Q_i$.  However, if the $\ell^d$-torsion points are not defined over $\F_{q}$, then the points $Q_i$ will span a proper subspace of $J[\ell^d]$.  If $\alpha$ is an endomorphism then $T$ will act as zero on all of the $Q_i$ and Algorithm \ref{a:frob} will output {\tt true}.   However, if $\alpha$ is not an endomorphism then $T$ may still act as zero on all of the $Q_i$ (in which case it must have nonzero action on the $\ell^d$-torsion points that are not defined over $\F_{q}$), and the algorithm will incorrectly output {\tt true}.  Thus to test whether $\alpha$ is an endomorphism, we must combine Algorithm \ref{a:frob} with a method of checking the field of definition of the $\ell^d$-torsion points, via the probabilistic method of Algorithm \ref{a:fields} or one of the other methods.
\end{rem}

\begin{prop}
\label{p:frobtime}
Let $J$ be the Jacobian of a genus 2 curve over $\F_p$.  
Assume that the zeta function of $J/\F_p$ is known, so that the cost to compute $\#J(\F_{p^k}) = \ell^s m$ is negligible.
Then the expected number of operations in $\F_p$ necessary to execute Algorithm \ref{a:frob} on $J/\F_{p^k}$ (ignoring $\log \log p$ factors) is 
\begin{equation*}
O(k^2  \log k (\log^2 p) \ell^{s - 4d}(-\log \epsilon))
\end{equation*}
\end{prop}

\begin{proof}
Let $q = p^k$.  In the proof of Proposition \ref{p:fieldtime}, we computed that the cost of computing a random point on $J(\F_q)$ is $O(\log q)$ operations in $\F_q$, and the cost of a point multiplication on $J(\F_q)$ is $O(\log q)$ operations in $\F_q$.  The chance that a random point in the $\ell$-primary part of $J(\F_{q})$ has exact order $\ell$ is $\frac{\ell^{4d} - \ell^{4d-4}}{\ell^s}$, so the expected number of random points necessary to find one point of exact order $\ell^d$ is $O(\ell^{s - 4d})$.  The cost of computing the Frobenius action is proportional to the cost of raising an element of $\F_q$ to the $p$th power, which is $O(\log p)$ $\F_q$-operations.

We conclude that the expected cost of a single trial with a random point is
\begin{equation*}
O(\log q + \log q + \log p) \ell^{s - 4d} M(q)
\end{equation*}
operations in $\F_p$, where $M(q)$ is the number of field operations in $\F_p$ needed to perform one field operation in $\F_q$.  If fast multiplication techniques are used, then $M(q) = O(\log q \log \log q) = O(k \log k \log p)$ (ignoring $\log \log p$ factors), so each trial takes 
$O(k^2  \log k (\log^2 p) \ell^{s - 4d})$
field operations in $\F_p$.  The number of points of exact order $\ell^d$ computed is $O(-\log \epsilon)$.  Putting this all together gives a total of
$O(k^2  \log k (\log^2 p) \ell^{s - 4d}(-\log \epsilon))$
field operations in $\F_p$.
\end{proof}

\subsection{The Couveignes method}
\label{ss:couveignes}

Recall that to test whether an element $\alpha \in \O_K$ of the form \eqref{eq:alpha-l} is an endomorphism of $J$, we determine whether the operator $T = a_0 + a_1 \pi + a_2 \pi^2 + a_3 \pi^3$ acts as zero on all elements of a set $\{Q_i\}$ that spans $J[\ell^d]$. 
Algorithm \ref{a:frob} computes the spanning 
set by choosing random points $P_i$ in $J(\F_{p^k})$, multiplying by 
an appropriate $m$ to get points $Q_i$ in the $\ell$-primary part of 
$J(\F_{p^k})$ (denoted $J(\F_{p^k})_\ell$), and keeping only those 
$Q_i$ whose order is exactly $\ell^d$.  If $J(\F_{p^k})_\ell$ is much 
larger than $J[\ell]$, the orders of most of the $Q_i$ will be too large, 
and it will take many trials to find the required number of points of 
order exactly $\ell^d$.  To reduce the number of trials required, we 
would like to find a function from $J(\F_{p^k})_\ell$ to $J[\ell^d]$ 
that sends most of the $Q_i$ to points of exact order $\ell^d$.

One way to compute such a function is as follows: compute the order 
$\ell^{t_i}$ of each $Q_i$; if $t_i \ge d$ send $Q_i \mapsto 
[\ell^{t-d}]Q_i$, otherwise send $Q_i \mapsto O$.  In most cases the 
image has order $\ell^d$.  However, since the multiplier 
$\ell^{t_i-d}$ will be different for each $Q_i$, this function does 
not define a group homomorphism, and thus the image of a set of 
points uniformly distributed in $J(\F_{p^k})_\ell$ will not be 
uniformly distributed in $J[\ell^d]$.

Couveignes \cite{couveignes} 
has described a map that has the properties we want and is a group 
homomorphism.  The idea is the following: if $\pi^k - 1 \in \ell^d 
\End(J)$, then there is an endomorphism $\phi$ such that $\ell^d \phi 
= \pi^k -1$.  Since $\pi^k-1$ acts as zero on $J(\F_{p^k})$, the 
image of $\phi$ on $J(\F_{p^k})$ must consist of $\ell^d$-torsion 
points.  Furthermore, the kernel of $\phi$ contains $\ell^d 
J(\F_{p^k})$, since $\phi(\ell^d P) = (\pi^k-1)(P) = 0$ if $P$ is 
defined over $\F_{p^k}$.  Thus we have a map
\begin{equation*}
\phi : 
J(\F_{p^k})/\ell^d J(\F_{p^k}) \to 
J[\ell^d].
\end{equation*}
Couveignes then uses the non-degeneracy of 
the Frey-R\"uck pairing (see \cite{schaefer}) to show that $\phi$ is 
a bijection.  Thus for any $Q_i$ not in $\ell J(\F_{p^k})$, 
$\phi(Q_i)$ has order exactly $\ell^d$.  Since $\phi$ is a surjective 
group homomorphism, the image of a set of points uniformly 
distributed in $J(\F_{p^k})$ will be uniformly distributed in 
$J[\ell^d]$.  The chance that $Q_i \in \ell J(\F_{p^k})$ is 
$1/\ell^{4}$, so applying $\phi$ to the $Q_i$ will very quickly give 
a spanning set of $J[\ell^d]$.

However, there is one important caveat: 
we may not be able to compute $\phi$.
The only endomorphisms we can 
compute are those involving the action of Frobenius and 
scalar multiplication; namely, endomorphisms in $\Z[\pi]$.  Thus we 
need to take $k$ to be the smallest integer such that $\pi^k - 1 \in 
\ell^d \Z[\pi]$.  We can then use the characteristic polynomial of 
Frobenius to write $\phi = \frac{\pi^k-1}{\ell^d} = M(\pi)$, where 
$M$ is a polynomial of degree $3$.  Furthermore, since we are applying 
$\phi$ only to points $Q_i \in J(\F_{p^k})_\ell$, we may reduce the 
coefficients of $M$ modulo $\ell^s$ and get the same action on the 
$Q_i$.

We have implemented the map $\phi$ in Magma and tested it on 
the examples that appear in Section \ref{s:ex}.  In our examples, the 
smallest $k$ for which $\pi^k - 1 \in \ell^d \Z[\pi]$ is usually 
equal to $\ell k_0$, where $k_0$ is the integer output by Algorithm 
\ref{a:torsdef}.  We found that the cost of choosing random points 
over a field of degree $\ell$ times as large far 
outweighs the benefit of having to reject fewer of the points $Q_i$, 
so this technique does not help to speed up Algorithm \ref{a:frob}.



\section{Bounding the field of definition of the $\ell^d$-torsion points}
\label{s:time}

The running times of Algorithms \ref{a:fields} and \ref{a:frob}  depend primarily on the size of the field $\F_{p^k}$ over which the $\ell^d$-torsion points of $J$ are defined.  In this section, we bound the size of $k$ in terms of $\ell^d$ and $p$.  We also show that to determine the field of definition of the $\ell^d$-torsion points of $J$ for $d > 1$, it suffices to determine the field of definition of the $\ell$-torsion points of $J$.  This result allows us to work over much smaller fields in Algorithm \ref{a:fields}, thus saving us a great deal of computation.

By Lemma \ref{l:dividesindex}, the prime powers $\ell^d$ input to Algorithms \ref{a:fields} and \ref{a:frob} divide the index $\Z[\pi,\pibar]$.  Thus a bound on this index gives a bound on the $\ell^d$ that appear.

\begin{prop}
\label{p:discbound}
Let $K$ be a primitive quartic CM field with discriminant $\Delta = \Delta(\O_K)$.  Suppose $\pi \in \O_K$ corresponds to the Frobenius endomorphism of the Jacobian of a genus 2 curve defined over $\F_p$.  Then
\begin{equation*}
[\O_K : \Z[\pi,\pibar]] \leq \frac{16 p^2}{\sqrt{\Delta}}.
\end{equation*}
\end{prop}

\begin{proof}
We showed in the proof of Corollary \ref{c:pibar} that $[\Z[\pi,\pibar]: \Z[\pi]] = p$.  Combining this result with the formula 
$$[\O_K : \Z[\pi]] = [\O_K : \Z[\pi,\pibar]] [\Z[\pi,\pibar] : \Z[\pi]], $$
we see that it suffices to show that $[\O_K : \Z[\pi]] \leq 16p^3/\sqrt{\Delta}$.  
(Note that $\Delta > 0$ by \cite[Proposition 9.4]{Howe}.)
Next, recall that $$[\O_K : \Z[\pi]] = \sqrt{\frac{\Disc(\Z[\pi])}{\Disc(\O_K)}}. $$
It thus suffices to show that $\sqrt{\Disc(\Z[\pi])} \leq 16 p^3$.  By definition,
\begin{equation}
\label{eq:disc}
 \sqrt{\Disc(\Z[\pi])} = \prod_{i < j} \abs{\alpha_i - \alpha_j},
\end{equation}
where $\alpha_i$ are the possible embeddings of $\pi$ into $\C$.  Since $\pi$ represents an action of Frobenius, all of the $\alpha_i$ lie on the circle $\abs{z} = \sqrt{p}$.  The product \eqref{eq:disc} takes its maximum value subject to this constraint when the $\alpha_i$ are equally spaced around the circle, which happens when the $\alpha_i$ are $\sqrt{p}$ times primitive eighth roots of unity.  The maximum product is thus $p^3 \sqrt{\Disc(\Q(\zeta_8))} = 16 p^3$.
\end{proof}

Proposition \ref{p:discbound} also follows directly from \cite[Proposition 7.4]{LPP}, where it is proved in a different manner that 
$\sqrt{\Disc(\Z[\pi, \pibar])} \leq 16 p^2$.

The next two propositions give tight bounds on the degree $k$ of the extension field of $\F_p$ over which the $\ell^d$-torsion points of $J$ are defined.  The first considers the case $d = 1$, and the second shows that as $d$ increases, $k$ grows by a factor of $\ell^{d-1}$.

\begin{prop}
\label{p:ltors}
Let $J$ be the Jacobian of a genus 2 curve over $\F_p$, and suppose that $\End(J)$ is isomorphic to the ring of integers $\O_K$ of the primitive quartic CM field $K$.  Let $\ell \neq p$ be a prime number, and suppose $\F_{p^k}$ is the smallest field over which the points of $J[\ell]$ are defined.  If $\ell$ is unramified in $K$, then $k$ divides one of the following:
\begin{itemize}
\item $\ell-1$, if $\ell$ splits completely in $K$;
\item $\ell^2-1$, if $\ell$ splits into two or three prime ideals in $K$;
\item $\ell^3 - \ell^2 + \ell - 1$, if $\ell$ is inert in $K$.
\end{itemize}
If $\ell$ ramifies in $K$, then $k$ divides one of the following:
\begin{itemize}
\item $\ell^3 - \ell^2$, if there is a prime over $\ell$ of ramification degree $3$, or if $\ell$ is totally ramified in $K$ and $\ell \leq 3$;
\item $\ell^2 - \ell$, in all other cases where $\ell$ factors into four prime ideals in $K$ (counting multiplicities);
\item $\ell^3 - \ell$, if $\ell$ factors into two or three prime ideals in $K$ (counting multiplicities).
\end{itemize}
\end{prop}

\begin{proof}
Let $\pi \in \O_K$ correspond to the Frobenius endomorphism.  By \cite[Fact 10]{el}, the $\ell$-torsion points of $J$ are defined over $\F_{p^k}$ if and only if $\pi^k-1 \in \ell\O_K$.  We observe that by the Chinese Remainder Theorem, this condition is satisfied if and only if $\pi^k \equiv 1 \pmod{\p_i^{e_i}}$ for all primes $\p_i \mid \ell\O_K$, where $e_i$ is the ramification degree of $\p_i$.  
Next, we note that the condition $\ell \neq p$ implies that $\pi \not \in \p_i$ for all $i$.  To see why this is true, suppose the contrary: $\pi \in \p_i$.  Since $\pi \pibar = p$, we have $p \in \p_i$, contradicting the fact that $\p_i$ is a prime over $\ell \neq p$.

From these observations we deduce that $k$ is the least common multiple of the multiplicative orders of $\pi$ mod each $\p_i^{e_i}$, and thus $k$ must divide the least common multiple of
$$\#(\O_K/\p_i^{e_i}\O_K)^\times = \ell^{f_i(e_i-1)}(\ell^{f_i}-1), $$
where $f_i$ is the inertia degree of $\p_i$.  We now consider the various possibilities for the splitting of $\ell$ in $\O_K$.  

First, suppose $\ell$ is unramified, so $e_i = 1$ for all $i$.
\begin{itemize}
\item If $\ell$ splits completely, then the inertia degrees of all the $\p_i$ are $1$, so $k \mid \ell - 1$.
\item If $\ell$ splits into two or three ideals, then at least one $\p_i$ has $f_i = 2$ and all have $f_i \leq 2$, so $k \mid \ell^2 - 1$.
\item If $\ell$ is inert, then there is a single $\p_i$ with $f_i = 4$, and $k$ divides $\ell^4 - 1$.  We will return to this case below to get a better bound.
\end{itemize}

Now suppose $\ell$ ramifies; there are six possibilities for the splitting of $\ell$ in $\O_K$.
\begin{itemize}
\item If $\ell \O_K = \p^3\q$, then $\p$ and $\q$ have inertia degree $1$, so $k$ divides $\ell^2(\ell-1)$.
\item If $\ell \O_K = \p^4$, then $\O_K/\p \cong \F_{\ell}$, and thus we have $\pi^{\ell-1} = 1 + \tau$ for some $\tau \in \p$.  There are now two subcases:
\begin{itemize}
	\item If $\ell \geq 5$, then $(1+\tau)^\ell \in 1 + \p^4$, so $\pi^{\ell(\ell-1)} \equiv 1 \pmod{\p^4}$.  Thus $k$ divides $\ell(\ell-1)$.  
	\item If $\ell = 2$ or $3$, then $(1+\tau)^\ell \equiv 1 + \tau^\ell \pmod{\p^4}$, so we must raise the expression to the $\ell$th power again to get rid of the $\tau^\ell$ term.  Thus $\pi^{\ell^2(\ell-1)} \equiv 1 \pmod{\p^4}$, and $k$ divides $\ell^2(\ell-1)$.
\end{itemize}
\item If $\ell \O_K = \p^2\q^2$ or $\p^2 \q \rr$, then all of the primes in question have inertia degree $1$, so $k$ divides $\ell(\ell-1)$.
\item If $\ell \O_K = \p^2\q$, then $\p$ has inertia degree $1$ and $\q$ has inertia degree $2$, so $k$ divides $\lcm( \ell(\ell-1), \ell^2-1) = \ell(\ell^2-1)$.
\item If $\ell \O_K = \p^2$, then $\O_K/\p \cong \F_{\ell^2}$, and thus we have $\pi^{\ell^2-1} = 1 + \tau$ for some $\tau \in \p$.  Then $(1+\tau)^\ell \in 1 + \p^2$, so $\pi^{\ell(\ell^2-1)} \equiv 1 \pmod{\p^2}$.  Thus $k$ divides $\ell(\ell^2-1)$.
\end{itemize}

Thus far we have used only the fact that $\pi$ is an algebraic integer, and we have not used the property that it represents the action of Frobenius.  To get a better bound in the case where $\ell$ is inert in $K$, we recall that since $\pi$ is the Frobenius endomorphism, we have $\pi\pibar = p$, and $K = \Q(\pi)$.  Since $\ell$ is inert, reduction modulo $\ell$ gives an injective group homomorphism 
$$\phi \colon \Aut\left(^{\displaystyle{K}}/_{\displaystyle{\Q}}\right) \to \Aut\left(^{\displaystyle{(\O_K/\ell\O_K)}}/_{\displaystyle{(\Z/\ell\Z)}}\right).$$  
Furthermore, the target group is isomorphic to $\Gal(\F_{\ell^4}/\F_\ell)$.  This group is cyclic of order $4$ and is generated by the $\ell$th-power Frobenius automorphism.  Since complex conjugation has order $2$ in $\Aut(K/\Q)$, its image under $\phi$ must be the map $\alpha \mapsto \alpha^{\ell^2}$. 
Thus $\pibar \equiv \pi^{\ell^2} \pmod{\ell}$, and $\pi^{\ell^2+1} \equiv p \pmod{\ell}$.  Since $p$ must reduce to an element of $\F_\ell^\times$, $p$ has order dividing $\ell-1$, so $\pi$ must have order dividing $(\ell^2+1)(\ell-1)$.
\end{proof}

The following proposition shows that in the cases we need for our application, the field of definition of the $\ell^d$-torsion points is determined completely by the field of definition of the $\ell$-torsion points.

\begin{prop}
\label{p:ldtors}
Let $A$ be an ordinary abelian variety defined over a finite field $F$, and let $\ell$ be a prime number not equal to the characteristic of $F$.  Let $d$ be a positive integer, and let $F'$ be the extension field of $F$ of degree $\ell^{d-1}$.  If the $\ell$-torsion points of $A$ are defined over $F$, then the $\ell^d$-torsion points of $A$ are defined over $F'$.  If $\End(A)$ is integrally closed, then the converse also holds.
\end{prop}

\begin{proof}
Let $R = \End(A)$, and let $\pi \in R$ be the Frobenius endomorphism of $F$.  By \cite[Fact 10]{el}, for any positive integers $t$ and $k$, the $\ell^t$-torsion points of $A$ are defined over the degree-$k$ extension of $F$ if and only if $\frac{\pi^k-1}{\ell^t} \in R$, i.e.\ $\pi^k \equiv 1 \pmod{\ell^t R}$.  To prove the proposition, it suffices to show that 
\begin{eqnarray*}
\pi \equiv 1 \pmod{\ell R} & \Leftrightarrow & \pi^{\ell^{d-1}} \equiv 1 \pmod{\ell^d R},
\end{eqnarray*}
with ($\Leftarrow$) holding when $R$ is integrally closed.

First suppose that $\pi^k \equiv 1 \pmod{\ell^t R}$, with $t \geq 1$.  Then we can write $\pi^k = 1 + \ell^t y$ for some $y \in R$.  Then 
\begin{equation*}
\pi^{k\ell} = 1 + \ell (\ell^t y) + \binom{\ell}{2} (\ell^t y)^2 + \cdots + (\ell^t y)^\ell,
\end{equation*}
so $\pi^{k\ell} \equiv 1 \pmod{\ell^{t+1} R}$.  We conclude that if the points of $A[\ell^t]$ are defined over the degree-$k$ extension of $F$, then the points of $A[\ell^{t+1}]$ are defined over the degree-$k\ell$ extension of $F$.  Thus if $A[\ell] \subset A(F)$, then by induction $A[\ell^d] \subset A(F')$.

Now suppose that $\pi^{k\ell} \equiv 1 \pmod{\ell^{t} R}$, with $t \geq 2$.  Since $A$ is ordinary, $R$ is an order in a number ring.  Thus if $R$ is integrally closed then it is a Dedekind domain, and we may write $\ell R = \prod \p_i^{e_i}$ uniquely for prime ideals $\p_i \subset R$.  By the Chinese Remainder Theorem, $\pi^k \equiv 1 \pmod{\ell^t R}$ if and only if $\pi^k \equiv 1 \pmod{\p_i^{e_i t}}$ for each $i$, so we may consider the problem locally at each $\p_i$.  Localizing and completing the ring $R$ at the prime $\p_i$ gives a complete local ring $R_v$ with maximal ideal $\p_i$ and valuation $v$ satisfying $v(\ell) = e_i$.

By hypothesis, we may write $\pi^{k \ell} = 1 + y$ for some $y \in \p_i^{e_i t}$.  We can define the $\ell$th-root function on $R_v$ to be
\begin{equation}
(1 + y)^{1/\ell} = \exp\left(\frac{1}{\ell} \log(1 + y) \right).
\end{equation}
By \cite[Proposition II.5.5]{neu}, if $y \in \p_i^{e_i t}$ then $\log(1 + y) \in \p_i^{e_i t}$.  Since $v(\ell) = e_i$, we have $v(\frac{1}{\ell}\log(1 + y)) \ge e_i(t-1)$, so by the same Proposition $(1 + y)^{1/\ell}$ converges and is in $1 + \p_i^{e_i(t-1)}$ whenever $(t-1)(\ell-1) > 1$.  Thus if $(t-1)(\ell-1) > 1$ then $\pi^k \equiv 1 \pmod{\p_i^{e_i(t-1)}}$.  We conclude that if $t > 2$ or $\ell > 2$ and the points of $A[\ell^t]$ are defined over the degree-$k\ell$ extension of $F$, then the points of $A[\ell^{t-1}]$ are defined over the degree-$k$ extension of $F$.  If $A[\ell^d] \subset A(F')$, then by descending induction $A[\ell] \subset A(F)$ if $\ell$ is odd, and $A[4] \subset A(F_2)$  if $\ell = 2$, where $F_2$ is the quadratic extension of $F$.

It remains to show that if $A[4] \subset A(F_2)$, then $A[2] \subset A(F)$.  This is equivalent to showing that if $\pi^2 - 1 \in 4R$ then $\pi-1 \in 2R$.  We prove the contrapositive: suppose $\pi-1 \not \in 2R$.  Then there is some prime $\p$ over $2$ such that $v_\p(\pi-1) < v_\p(2)$.  Since $\pi+1 = (\pi-1) +2$ and $v_\p(\pi-1) < v_\p(2)$, we must also have $v_\p(\pi+1) < v_\p(2)$.  Multiplying the two expressions gives $v_\p(\pi^2-1) < v_\p(4)$, so $\pi^2 - 1$ cannot be contained in $4R$.  We conclude that $\pi^2 -1 \in 4R$ implies $\pi-1 \in 2R$.
\end{proof}

\begin{cor}
\label{c:ldtors}
Let $J$ be the Jacobian of a genus 2 curve over $\F_p$, and suppose that $\End(J)$ is isomorphic to the ring of integers $\O_K$ of the primitive quartic CM field $K$.   Let $\ell^d$ be a prime power with $\ell \neq p$, and suppose $\F_{p^k}$ is the smallest field over which the points of $J[\ell^d]$ are defined.  Then $k < 3 p^6$.
\end{cor}

\begin{proof}
By Proposition \ref{p:ltors}, the points of $J[\ell]$ are defined over a field $F$ of degree less than $\ell^3$ over $\F_p$.  By Proposition \ref{p:ldtors}, the points of $J[\ell^d]$ are defined over a field $L$ of degree $\ell^{d-1}$ over $F$.  Since degrees of extensions multiply, we get 
\[ k = [L : \F_p] < \ell^{d+2} \leq \ell^{3d}. \]
By Proposition \ref{p:discbound}, $\ell^d \le \frac{16}{\sqrt{\Delta}} p^2$, where $\Delta$ is the discriminant of the quartic CM field $K$.  
Lemma \ref{l:cmdisc} below shows that any primitive quartic CM field has $\Delta \ge 125$,
so $\ell^d \leq \frac{16}{\sqrt{125}} p^2$.
Since $k < \ell^{3d}$, we conclude that $k < 3 p^6$.
\end{proof}


\begin{lemma}
\label{l:cmdisc}
Suppose $K$ is a primitive quartic CM field.  Then $\Delta(K) \ge 125$.
\end{lemma}

\begin{proof}
Since $\Delta(\Q(\zeta_5)) = 125$, it suffices to show that no smaller discriminant can occur.  The fact that $\Disc(K) > 0$ follows from \cite[Proposition 9.4]{Howe}.  Now suppose $\Delta(K) < 125$.  Since $\Delta(K_0)^2 \mid \Delta(K)$, we must have $K_0 = \Q(\sqrt{2})$ or $\Q(\sqrt{5})$, as these are the only two real quadratic fields with discriminant less than $12$.  Since $\Q(\sqrt{2})$ has class number $1$, by \cite[Proposition VI.6.9]{neu}, $\Q(\sqrt{2})$ has no unramified quadratic extensions, so 
$\Delta(K)$ is strictly greater than $\Delta(K_0)^2$.  Thus if $K_0 = \Q(\sqrt{2})$ then $\Delta(K) \ge 128$.  

We deduce that $K_0 = \Q(\sqrt{5})$ and $K$ must be of the form $\Q(i\sqrt{a+b\sqrt{5}})$, with $a$, $b$, and $a^2 - 5b^2$ positive integers.  
Since $K$ is primitive, $a^2 - 5b^2$ is not a square in $\Q$ and its square-free part divides $\Delta(K)/\Delta(K_0)^2$.  It thus suffices to show that the square-free part of $a^2 - 5b^2$ is at least $5$; this follows from the fact that $2$ and $3$ are inert in $\Q(\sqrt{5})$, so there are no integer solutions to $a^2 - 5b^2 = 2$ or $3$.  
\end{proof}




\section{Computing Igusa class polynomials}
\label{s:crt}

This section combines the results of all of the previous sections into a full-fledged probabilistic version of Eisentr\"ager and Lauter's CRT algorithm to compute Igusa class polynomials for primitive quartic CM fields \cite[Theorem 1]{el}.  

\begin{alg}
\label{a:crt}
The following algorithm takes as input a primitive quartic CM field $K$, three integers $\lambda_1,\lambda_2,\lambda_3$ which are multiples of the denominators of the three Igusa class polynomials, and a real number $\epsilon > 0$, and outputs three polynomials $H_1,H_2,H_3 \in \Q[x]$.  With high probability, the polynomials $H_i(x)$ output by the algorithm are the Igusa class polynomials for $K$.
\begin{enumerate} 
\item{(Initialization.)}

\begin{enumerate}
\item Let $D$ be the degree of the Igusa class polynomials for $K$, computed via class number algorithms, e.g.
\cite[Algorithm 6.5.9]{coh}.
\item Compute an integral basis $\B$ for $\O_K$, using e.g.~\cite[Algorithm 6.1.8]{coh}.
\item Set $p \leftarrow 3$, $B \leftarrow 1$, $H_1,H_2,H_3 \leftarrow 0$, $F_1,F_2,F_3 \leftarrow 0$.
\end{enumerate}

\item \label{st:increase} Set $p \leftarrow \mathtt{NextPrime}(p)$ until $p$ splits completely in $K$ and $p$ splits into principal ideals in $K^*$ (the reflex field of $K$).  

\item \label{st:igusaloop} (Finding the curves.) Set $T_1,T_2,T_3 \leftarrow \{\}$.  For each $(i_1,i_2,i_3) \in \F_p^3$, do the following:
	\begin{enumerate}
	\item \label{st:computecurve} Compute a curve $C/\F_p$ with Igusa invariants $(i_1,i_2,i_3)$, using the algorithms of Mestre \cite{mes} and Cardona-Quer \cite{cq}.
	\item \label{st:zeta} 
	Run Algorithm \ref{a:zeta} with inputs $K, p, C$.  
		\begin{enumerate}
		\item If the algorithm outputs {\tt false}, go to the next triple $(i_1,i_2,i_3)$.  
		\item If the algorithm outputs {\tt true}, let $\pi$ be one of the possible Frobenius elements it outputs.  
		\end{enumerate}
	\item \label{st:computefields} 
	For each prime $\ell$ dividing $[\O_K : \Z[\pi,\pibar]]$, do the following:
		\begin{enumerate}
		\item \label{st:torsdef} Run Algorithm \ref{a:torsdef} with inputs $K$, $\ell$, $\pi$.  Let the output be $k$.
		\item \label{st:torsdef2} Run Algorithm \ref{a:fields} with inputs $\Jac(C)$, $\F_{p^k}$, $\ell$, and $\epsilon$.  If the output is {\tt false}, go to the next triple $(i_1,i_2,i_3)$.
		\item \label{st:computefrob} If $\ell^2$ divides $[\O_K : \Z[\pi,\pibar]]$, then for each $\alpha \in \B \setminus \Z$ written in the form \eqref{eq:alpha} with denominator $n$, do the following:
			\begin{enumerate}
			\item Let $d$ be the largest integer such that ${\ell^d} \mid n$.  If $d=0$, go to the next $\alpha$. 
			\item Set $k' \leftarrow k \ell^{d-1}$.
			\item Run Algorithm \ref{a:frob} with inputs $\Jac(C), \F_{p^{k'}}, \ell^d, \pi, \frac{n}{\ell^d}\alpha, \epsilon$.  
			\item If Algorithm \ref{a:frob} outputs {\tt false}, go to the next triple $(i_1,i_2,i_3)$.  Otherwise go to the next $\alpha$.
			\end{enumerate}
		\end{enumerate}	
	\item Adjoin $i_1,i_2,i_3$ to the sets $T_1,T_2,T_3$, respectively (counting multiplicities).
	\end{enumerate}
\item \label{st:toomany} If the size of each set $T_1,T_2,T_3$ is not equal to $D$, go to Step \ref{st:increase}.
\item \label{st:polynomials} (Computing the Igusa class polynomials.) For $i \in \{1,2,3\}$, do the following:
	\begin{enumerate}
	\item \label{st:combine} Compute $F_{i,p}(x) = \lambda_i \prod_{j \in T_i} (x - j)$ in $\F_p[x]$.
	\item \label{st:crtcombine} Use the Chinese Remainder Theorem to compute $F_i'(x) \in \Z[x]$ such that $F_i'(x) \equiv F_i(x) \pmod{B}$, $F_i'(x) \equiv F_{i,p}(x) \pmod{p}$, and the coefficients of $F_i'(x)$ are in the interval $[-pB/2,pB/2]$.
	\item \label{st:output} If $F_i'(x) = F_i(x)$, output $H_i(x) = \lambda_i^{-1} F_i(x)$.  If $H_i(x)$ has been output for all $i$, terminate the algorithm.
	\item Set $F_i(x) \leftarrow F_i'(x)$.
	\end{enumerate}
\item Set $B \leftarrow pB$, and return to Step \ref{st:increase}.
\end{enumerate}

\end{alg}

\begin{proof}
In view of \cite[Theorem 1]{el}, it suffices to prove that Step \ref{st:computefields} correctly determines the set of curves with $\End(\Jac(C)) = \O_K$.  It follows from Section \ref{s:generators} that $\End(\Jac(C)) = \O_K$ if and only if each of the elements of the generating set listed in Proposition \ref{prop:generators} is an endomorphism.  

By Algorithm \ref{a:zeta}, the $\pi$ computed in Step \ref{st:zeta} is such that $\pi^\sigma$ is the Frobenius element of $\Jac(C)$ for some $\sigma \in \Aut(K/\Q)$.  By Corollary \ref{c:galois-generators}, $\End(\Jac(C)) = \O_K$ if and only if $\beta^\sigma$ is an endomorphism for each $\beta$ in the generating set of Proposition \ref{prop:generators}.  
Since elements of $\Aut(K/\Q)$ preserve $\O_K$ as a set, $[\O_K : \Z[\pi^\sigma,\pibar^\sigma]] = [\O_K : \Z[\pi,\pibar]]$.

For each $\ell$ dividing $[\O_K : \Z[\pi,\pibar]]$, Steps \ref{st:torsdef} and \ref{st:torsdef2} test probabilistically whether $\frac{(\pi^\sigma)^k -1}{\ell}$ is an endomorphism for an appropriate $k$.  By Corollary \ref{c:nofrob}, for any such $\ell$ dividing $[\O_K : \Z[\pi,\pibar]]$ exactly, this suffices to determine whether $\frac{n}{\ell}\alpha^\sigma$ is an endomorphism for each $\alpha \in \B \setminus \Z$.  

By Corollary \ref{c:galois-frob}, if $\ell^2$ divides $[\O_K : \Z[\pi,\pibar]]$ then Step \ref{st:computefrob} tests probabilistically whether $\frac{n}{\ell^d}\alpha^\sigma$ is an endomorphism.  The input uses the field $\F_{p^{k'}}$ because Proposition \ref{p:ldtors} implies that if the $\ell$-torsion points are defined over $\F_{p^k}$, then the $\ell^d$-torsion points are defined over $\F_{p^{k'}}$.
\end{proof}

\begin{rem}
\label{r:crtstep}  Note that Step \ref{st:polynomials} differs from the corresponding step in \cite[Theorem 1]{el}.
Our version of the algorithm minimizes the amount of computation by terminating the algorithm in Step \ref{st:output} as soon as the polynomials agree modulo two consecutive primes. 
For each prime $p_i$ in an increasing sequence of primes, we compute a polynomial $F_{i,p}(x)$ that is congruent to the Igusa class polynomial $H_i(x)$ modulo the prime $p_i$ \cite[Theorem 2]{el}. 
We then use the Chinese Remainder Theorem and the collection of polynomials $\{F_{i,p}(x) \}$ to compute a polynomial $F_i(x)$ modulo $b_i = \prod_{j = 1}^{i} p_j$.  
If $F_i(x)=F_{i+1}(x)$, then with high probability the coefficients of $H_i(x)$ are less than $b_{i+1}$, and thus $F_{i}(x)$ is equal to $H_i(x)$ itself.
This conclusion is justified by the fact that if an integer $n$ has the property that it is the same modulo $b_i$ and 
modulo $b_{i+1}$, then $n= a_i + r_i b_i = a_{i+1} + r_{i+1}b_{i+1}$, with $a_i < b_i$ and $a_i=a_{i+1}$.  It follows that $p_{i+1}$ divides $r_i$.  Since the probability of this happening for a random number $r_i$ is $1/p_{i+1}$, the probability
that all coefficients would simultaneously satisfy this congruence is $(1/p_{i+1})^{D+1}$, so most likely we have that actually $r_{i+1}=0$ for each coefficient.
\end{rem}

\begin{rem}
The $\lambda_i$ input into the algorithm can be taken to be products of primes bounded in \cite{goren-lauter}, raised to a power that will be made explicit in forthcoming work.  In practice, the power can be taken to be a small multiple of $6$.
\end{rem}  

Since we check after every prime $p_i$ whether the algorithm is finished, we do not need to know in advance the number of primes $p_i$ that we will need to use.  Thus the only bounds that need to be computed in advance are the bounds $\lambda_i$ on the denominators of the coefficients of the Igusa class polynomials.  In particular, we do not need to have a bound on either the numerators or the absolute values of the coefficients.


\section{Implementation notes}
\label{s:notes}

Our most significant observation is that in practice,  the running time of the probabilistic CRT algorithm is dominated by generating $p^3$ curves for each small $p$.  Steps \eqref{st:computecurve} and \eqref{st:zeta} of Algorithm \ref{a:crt} generate a list of curves $C$ for which $\End(\Jac(C))$ is an order in $\O_K$.  Algorithms \ref{a:fields} and \ref{a:frob} determine which endomorphism rings are equal to $\O_K$.  Data comparing the relative speeds of these two parts of the algorithm appear in Section \ref{s:ex}.  This section describes a number of ways to speed up Algorithm \ref{a:crt}, which are reflected in the running times that appear in Section \ref{s:ex}.

\begin{enumerate}
\item \label{n:skip} If $p$ and $k$ are large, then arithmetic on $J(\F_{p^k})$ is prohibitively slow, which slows down  Algorithms \ref{a:fields} and \ref{a:frob}.  Since for various $\ell$ dividing the index $[\O_K:\Z[\pi,\pibar]]$, the extension degrees $k$ depend only on the prime $p$ and the CM field $K$ and not on the curve $C$, these extension degrees may be computed in advance (via Algorithm \ref{a:torsdef}) before generating any curves.  We set some bound $N$ and tell the program that if the extension degree $k$ for some $\ell$ is such that $p^k > N$, we should skip that $p$ and go on to the next prime.  For example, if $K = \Q(i\sqrt{13+2\sqrt{13}})$ and $p = 53$ (see Example \ref{ex:13.2.13}), we have $[\O_K: \Z[\pi,\pibar]] = 3^2 \cdot 43$, and the $43$-torsion of a Jacobian $J$ with $\End(J) = \O_K$ will be defined over $\F_{p^{924}}$, a field of over 5000 bits that is far too large for our current implementation to handle efficiently.
\item In a similar vein, since the speed of Algorithms \ref{a:fields} and \ref{a:frob} is determined by the size of the fields $\F_{p^k}$, for optimum performance one should perform these calculations in order of increasing $k$,
so that as the fields get larger there are fewer curves to check.
\item  Algorithms \ref{a:fields} and \ref{a:frob} take a single curve as input.  In Algorithm \ref{a:crt} those algorithms are executed with the same field $K$ and many different curves, so any parameter that only depends on the field $K$ and the prime $p$ can be precomputed and stored for repeated reference.  For example, the representation $\alpha = (a_0 + a_1 \pi + a_2 \pi^2 + a_3 \pi^3)/n$ and the extension degrees $k$ in Step \ref{st:torsdef} can be computed only once.  In addition, all of the curves that pass Step \ref{st:zeta} have one of a small number of given zeta functions.  Since $\#J(\F_{p^k})$ is determined by the zeta function, this number can also be computed in advance.
\item If $\F_{p^k}$ is small enough, it may be faster to check fields of definition using the brute force method of Section \ref{ss:brute1}, rather than Algorithm \ref{a:fields}.  If $\ell$ is small (as must be the case for $k$ to be small), then we often find that $\#J(\F_{p^k}) = \ell^s m$ with $s \gg 4d$, and thus the number of random points needed in Algorithms \ref{a:fields} and \ref{a:frob} will be very large.  While computing the group structure is an exponential-time computation, we find that if the group has size at most $2^{200}$, MAGMA can compute the group structure fairly quickly.
\item \label{note:output} If Step \ref{st:output} has already output $H_j(x)$ for some $j$,  the roots of this polynomial mod $p$ can be used as the possible values of $i_j$ in Step \ref{st:igusaloop}.  This will greatly speed up the calculation of the $F_{i,p}$ for the remaining primes: if one $H_j$ has been output then only $p^2D$ curves need to be computed (instead of $p^3$), and if two $H_j$ have been output then only $pD^2$ curves need to be computed.
\item In practice, for small primes $p$ ($p < 800$ in our MAGMA implementation), computing $\#C(\F_p)$ (Step \ref{st:countc} of Algorithm \ref{a:zeta}) is more efficient than choosing a random point on $J(\F_p)$ and determining whether it is killed by one of the potential group orders (Step \ref{st:randompt} of Algorithm \ref{a:zeta}), so these two steps should be switched for maximum speed.  However, as $p$ grows, the order of the steps as presented will be the fastest.
\end{enumerate}


\section{Examples}
\label{s:ex}

This section describes the performance of Algorithm \ref{a:crt} on three quartic CM fields: $\Q(i\sqrt{2 + \sqrt{2}})$, $\Q(i\sqrt{13 + 2\sqrt{13}})$, and $\Q(i\sqrt{29 + 2\sqrt{29}})$.  These fields are all Galois and have class number $1$, so the density of primes with the desired splitting behavior is maximal.  The Igusa polynomials are linear; they have integral coefficients for the first two fields, and have denominators dividing $5^{12}$ for the last.  
In all three examples, as $p$ grows, the running time of the algorithm becomes dominated by the computation of $p^3$ curves for each $p$, whereas it was previously suspected that the the endomorphism ring computation would be the slow step in the CRT algorithm.  A fast implementation in C to produce the curves from their Igusa invariants and to test the numbers of points would thus significantly improve the running time of the CRT algorithm.

Details of the algorithms' execution are given below.  The algorithms were run on a 2.39 GHz AMD Opteron with 4 GB of RAM.  The table headings have the following meaning:
\begin{itemize}
\item $p$: Size of prime field over which curves were generated.
\item $\ell^d$: Prime powers appearing in the denominators $n$ of elements $\alpha$ input into Algorithms \ref{a:fields} and \ref{a:frob}, when written in the form \eqref{eq:alpha}.
\item $k$: Degrees of extension fields over which $\ell^d$-torsion points are expected to be defined.  These are listed in the same order as the corresponding $\ell^d$.
\item Curves: Time taken to generate $p^3$ curves and determine which have CM by $K$ (cf.\ Algorithm \ref{a:zeta}).
\item \#Curves: Number of curves computed whose Jacobians have CM by $K$.
\item \ref{a:fields} \& \ref{a:frob}: Time taken to run Algorithms \ref{a:fields} and \ref{a:frob} to find the single curve whose Jacobian has endomorphism ring equal to $\O_K$.
\end{itemize}

\begin{ex} 
\label{ex:2.1.2} 
We ran Algorithm \ref{a:crt} with $K = \Q(i\sqrt{2 + \sqrt{2}})$ and $\lambda_1,\lambda_2,\lambda_3 = 1$.  The results appear in Table 1.  The last column of the table shows the intermediate polynomials $F_i(x)$ computed via the Chinese Remainder Theorem in Step \ref{st:crtcombine}.  The algorithm output the $F_i(x)$ listed for $p = 151$ as the Igusa class polynomials of $K$.
\end{ex}

\begin{table}[ht]
\caption{Results for Algorithm \ref{a:crt} run with $K = \Q(i\sqrt{2 + \sqrt{2}})$ and $\lambda_1,\lambda_2,\lambda_3 = 1$.}
\begin{tabular}{|c|c|c|r|c|r|c|}
\hline
$p$ & $\ell^d$ & $k$ & Curves & \#Curves & \ref{a:fields} \& \ref{a:frob} & $F_i(x)$ \\
\hline \hline
7 	& 2,4 	& 2,4 		& 0.5 sec	& 7		& 0.3 sec  & \scriptsize
$\begin{array}{c} x + 2 \\ x + 5 \\ x + 6 \\ 
\pmod{7} \end{array}$ \\ \hline
17 	& 4,8 	& 2,4		& 4 sec 		& 39 	& 0.2 sec & \scriptsize
$\begin{array}{c} x - 54 \\ x + 19 \\ x - 8 \\ 
\pmod{119} \end{array}$ \\ \hline
23 	& 2,4,7	& 2,4,3		& 9 sec 		& 49 	& 2.3 sec  & \scriptsize
$\begin{array}{c} x + 1017 \\ x + 852 \\ x + 111 \\ 
\pmod{2737} \end{array}$ \\ \hline
71 	& 2,4 	& 2,4		& 255 sec 	& 7 		& 0.7 sec  & \scriptsize
$\begin{array}{c} x - 75619 \\ x + 28222 \\ x - 46418 \\ 
\pmod{194327} \end{array}$ \\ \hline
97 	& 4,8	& 2,4		& 680 sec 	& 39 	& 0.3 sec & \scriptsize
$\begin{array}{c} x - 8237353 \\ x + 9355918 \\ x + 9086951 \\ 
\pmod{18849719} \end{array}$ \\ \hline
103 & 2,4,17 & 2,4,16	& 829 sec 	& 119 	& 17.6 sec  & \scriptsize
$\begin{array}{c} x + 104860961 \\ x - 28343520 \\ x - 9762768 \\ 
\pmod{1941521057} \end{array}$\\ \hline
113 & 7,8,32	 & 6,4,16	& 1334 sec 	& 1281 	& 28.8 sec  & \scriptsize
$\begin{array}{c} x - 1836660096 \\ x - 28343520 \\ x - 9762768 \\ 
\pmod{219391879441} \end{array} $\\ \hline
151 & 2,4,7,17 & 2,4,6,16	& 0.2 sec & 1	& --  		& \scriptsize
$\begin{array}{c}  x - 1836660096 \\ x - 28343520 \\ x - 9762768 \\ 
\pmod{33128173795591} \end{array}$ \\ \hline
\end{tabular}
\end{table}

The total time of this run was 3162 seconds, or about 53 minutes.  We observe that the polynomials $F_2$ and $F_3$ agree for $p = 103$ and $p = 113$.  We deduce that these polynomials are the correct Igusa polynomials, and following note \eqref{note:output} of Section \ref{s:notes}, we use their roots for the values of $i_2$ and $i_3$ for $p = 151$.  Thus instead of computing $151^3 \approx 2^{22}$ curves, we need to compute only $151$ curves, out of which we can easily choose the right one.  As a result, the computation for $p = 151$ takes practically no time at all.  The same phenomenon also appears for the last prime in Examples \ref{ex:13.2.13} and \ref {ex:29.2.29}.

\begin{ex} 
\label{ex:13.2.13} 
We ran Algorithm \ref{a:crt} with $K = \Q(i\sqrt{13 + 2\sqrt{13}})$ and $\lambda_1,\lambda_2,\lambda_3 = 1$.  The results appear in Table 2.  The algorithm output the following Igusa class polynomials:
\[ x - 1836660096, \quad x - 28343520 \quad x - 9762768. \]

The total time of this run was 6969 seconds, or about 116 minutes.  In this example we skip some primes because Algorithms \ref{a:fields} and \ref{a:frob} would need to compute in fields which are too large to be practical.  In particular, for $p = 29,53,107,139$, the algorithms would run over extension fields of degree $264,924,308,162$, all of which have well over $1000$ bits.  Skipping these primes has no effect on the ultimate outcome of the algorithm.
\end{ex}

\begin{table}[ht]
\caption{Results for Algorithm \ref{a:crt} with $K = \Q(i\sqrt{13 + 2\sqrt{13}})$ and $\lambda_1,\lambda_2,\lambda_3 = 1$.}
\begin{tabular}{|c|c|c|r|c|r|}
\hline
$p$ & $\ell^d$ & $k$ & Curves & \#Curves & \ref{a:fields} \& \ref{a:frob} \\
\hline \hline
29 	& 3,23	& 2,264	& -- 		& -- 	& -- 		 \\ \hline
53 	& 3,43	& 2,924	& -- 		& -- 	& -- 		 \\ \hline
61 	& 3		& 2 		& 167 sec 	& 9 		& 0.2 sec 	 \\ \hline
79 	& 27		& 18 	& 376 sec 	& 81 	& 8.1 sec  	\\ \hline
107 	& 9,43	& 6,308	& -- 		& -- 	& -- 		\\ \hline
113 	& 3,53	& 1,52	& 1118 sec 	& 159 	& 137.2 sec  \\ \hline
131 & 9,53	& 6,52	& 1872 sec 	& 477 	& 127.4 sec   \\ \hline
139 & 9,243	& 6,162 	& -- 		& -- 	& -- 		 \\ \hline
157 & 9,81	& 6,54	& 3147 sec 	& 243 	& 16.5 sec 	 \\ \hline
191 & 3,4,8	& 2,2,4	& 0.2 sec 	& 1 		& -- 		\\ \hline
\end{tabular}
\end{table}

\begin{ex} 
\label{ex:29.2.29} 
We ran Algorithm \ref{a:crt} with $K = \Q(i\sqrt{29 + 2\sqrt{29}})$ and $\lambda_1,\lambda_2,\lambda_3 = 5^{12}$.  The results appear in Table 3.  The algorithm output the following Igusa class polynomials:
\[ x - {\textstyle \frac{2614061544410821165056}{5^{12}}}, \quad
 x + {\textstyle \frac{586040972673024}{5^6}}, \quad 
 x + {\textstyle \frac{203047103102976}{5^6}}.
\]

The total time of this run was 56585 seconds, or about 15 hours, 43 minutes.  In this example we again skip some primes because the fields input to Algorithms \ref{a:fields} and \ref{a:frob} would be too large.  We also note that for $p = 7$, $\O_K = \Z[\pi,\pibar]$, so any curve over $\F_7$ that has a correct zeta function already has CM by all of $\O_K$, and we do not need to run Algorithms \ref{a:fields} and \ref{a:frob}.

\end{ex}

\begin{table}[htb]
\caption{Results for Algorithm \ref{a:crt} with $K = \Q(i\sqrt{29 + 2\sqrt{29}})$ and $\lambda_1,\lambda_2,\lambda_3 = 5^{12}$.}
\begin{tabular}{|c|c|c|r|c|r|}
\hline
$p$ & $\ell^d$ & $k$ & Curves & \#Curves & \ref{a:fields} \& \ref{a:frob} \\
\hline \hline
7 & -- & -- & 0.3 sec & 1 & -- 
\\ \hline
23 & 13 & 84 & 9 sec & 15 & 70.7 sec 
\\ \hline
53 & 7 & 6 & 105 sec & 7 & 0.5 sec 
\\ \hline
59 & 4,5,8	& 2,12,4	& 164 sec & 322 & 6.4 sec 
\\ \hline
83 & 3,5	& 4,24	& 431 sec & 77 & 9.8 sec  
\\ \hline
103 & 67 & 1122 & -- & -- & --  \\ \hline
107 & 7,13 & 6,42 & 963 sec & 105 & 69.3 sec 
\\ \hline
139 & 7,25 & 2,60 & 2189 sec & 259 & 62.1 sec 
\\ \hline
181 & 9,27 & 6,18  & 84 min & 161 & 3.6 sec 
\\ \hline
197 & 5,109 & 24,5940 & -- & -- & -- 
\\ \hline
199 & 25 & 60 & 106 min & 37 & 1355.3 sec 
\\ \hline
223 & 4,8,23 & 2,4,22 & 174 min & 1058 & 35.1 sec 
\\ \hline
227 & 109 & 1485 & -- & -- & --  \\ \hline
233 & 5,7,13 & 8,3,28 & 193 min & 735 & 141.6 sec 
\\ \hline
239 & 7,109 & 6,297 & -- & -- & --  \\ \hline
257 & 3,7,13 & 4,6,84 & 286 min & 1155 & 382.8 sec 
\\ \hline
277 & 5,7,23 & 24,6,22 & 0.3 sec & 1 & -- 
\\ \hline
\end{tabular}
\end{table}


\newpage

\begin{rem}
The data in Examples \ref{ex:2.1.2}, \ref{ex:13.2.13}, and \ref{ex:29.2.29} suggest that odd primes dividing the index $[\O_K : \Z[\pi,\pibar]]$
always split in $\O_{K_0}$, the ring of integers of $K_0$. In fact the factorization of the index $[\O_K : \Z[\pi,\pibar]]$ was given in \cite[Proposition 5]{el} for primitive quartic CM fields $K$ when $K_0$ has class number 1.
We write $\pi = c_1 + c_2 \sqrt{d} + (c_3 + c_4\sqrt{d})\eta$, where the $c_i$ are rational numbers with only powers of $2$ in the denominators and $\eta = i\sqrt{a+b\sqrt{d}}$ with $a,b,d \in \Z$, $d>0$ and square-free.
Then the index is, up to powers of $2$, the product of $c_2$ with $(c_3^2 -c_4^2d)$, where $c_2$ is the index of $\Z[\pi+\pibar]$ in $\O_{K_0}$ up to a power of $2$.  If a prime divides $(c_3^2 -c_4^2d)$ exactly, i.e. the square of the prime does not divide it, then the prime splits in $K_0$.
Thus primes different from $2$ dividing the index $[\O_K : \Z[\pi,\pibar]]$ exactly either split in $K_0$ or divide the index $[\O_{K_0}: \Z[\pi+\pibar]$. So except possibly for primes dividing $c_2$, no odd primes dividing the index $[\O_K : \Z[\pi,\pibar]]$ exactly are inert or totally ramified in $K$.
If $K$ is Galois, then this is enough to ensure that the extension degree $k$ determined by Proposition \ref{p:ltors} is at most $\ell^2$. This agrees with the data in our examples, all of which considered Galois fields.

In practice, if a prime $\ell$ is inert or totally ramified in $K$, it would almost certainly be skipped anyway, since Proposition \ref{p:ltors} shows that the $\ell$-torsion may be defined over an extension field of degree $k \sim \ell^3$, which is too large to be practical (cf.~Note \eqref{n:skip} of Section \ref{s:notes}).
However the theoretical running times of Algorithms \ref{a:fields} and \ref{a:frob}, given by Propositions \ref{p:fieldtime} and \ref{p:frobtime} respectively, improve if inert or ramified primes $\ell$ are not considered.  The slow step of both algorithms is computing a random point on $J(\F_{p^k})$, which takes roughly $O(k^2 \log k (\log p)^2)$ $\F_p$ operations.  Since the bound on $\ell$ is $p^2$, if $k$ is bounded by $\ell^2$ instead of $\ell^3$, this step would run in $O(p^{8} \log^3 p)$ instead of $O(p^{12} \log^3 p)$ time.

\end{rem}

\end{document}